\documentclass{IEEEtran}

\usepackage{enumitem,kantlipsum}
\usepackage{graphicx}
\usepackage{epstopdf}
\usepackage{amssymb}
\usepackage{amsmath}
\usepackage{subfigure}
\usepackage{epsfig}
\usepackage{color}
\usepackage{enumitem}
\usepackage{float}
\usepackage{soul}
\usepackage{wrapfig}
\usepackage{arydshln}

\usepackage{amsthm}

\usepackage{url}

\newtheorem{Remark 1}{Remark}
\newtheorem{Remark 2}[Remark 1]{Remark}
\newtheorem{Remark 3}[Remark 1]{Remark}
\newtheorem{Remark 4}[Remark 1]{Remark}
\newtheorem{Remark 5}[Remark 1]{Remark}
\newtheorem{Remark 6}[Remark 1]{Remark}
\newtheorem{Remark 7}[Remark 1]{Remark}

\newtheorem{Lemma 1}{Lemma}
\newtheorem{Lemma 2}[Lemma 1]{Lemma}
\newtheorem{Lemma 3}[Lemma 1]{Lemma}
\newtheorem{Lemma 4}[Lemma 1]{Lemma}
\newtheorem{Lemma 5}[Lemma 1]{Lemma}
\newtheorem{Lemma 6}[Lemma 1]{Lemma}
\newtheorem{Lemma 7}[Lemma 1]{Lemma}

\newtheorem{Corollary 1}{Corollary}
\newtheorem{Corollary 2}[Corollary 1]{Corollary}

\newtheorem{Assumption 1}{Assumption}
\newtheorem{Assumption 2}[Assumption 1]{Assumption}
\newtheorem{Assumption 3}[Assumption 1]{Assumption}
\newtheorem{Assumption 4}[Assumption 1]{Assumption}
\newtheorem{Definition 1}{Definition}
\newtheorem{Theorem 1}{Theorem}
\newtheorem{Theorem 2}[Theorem 1]{Theorem}
\newtheorem{Theorem 3}[Theorem 1]{Theorem}
\newtheorem{Theorem 4}[Theorem 1]{Theorem}
\newtheorem{Theorem 5}[Theorem 1]{Theorem}
\newtheorem{Theorem 6}[Theorem 1]{Theorem}
\newtheorem{Theorem 7}[Theorem 1]{Theorem}
\newtheorem{Theorem 8}[Theorem 1]{Theorem}
\newtheorem{Theorem 9}[Theorem 1]{Theorem}
\newtheorem{Theorem 10}[Theorem 1]{Theorem}

\newtheorem{Proposition 1}{Proposition}
\newtheorem{Proposition 2}[Proposition 1]{Proposition}

\title{\LARGE \bf
 Quantization enabled  Privacy Protection in Decentralized Stochastic Optimization}

\author{Yongqiang Wang, Tamer Ba{\c{s}}ar
\thanks{ The work of the first author  was supported in part by the National Science Foundation under Grants  ECCS-1912702, CCF-2106293 and CCF-CCF-2215088.  Research of the second author was supported in part by the ONR MURI Grant N00014-16-1-2710 and in part by the Army Research Laboratory, United States, under Cooperative Agreement Number W911NF-17-2-0196.}
\thanks{Yongqiang Wang is with the Department of Electrical and Computer Engineering, Clemson University, Clemson, SC 29634, USA
{\tt\small{yongqiw}@clemson.edu}
}%
\thanks{Tamer Ba{\c{s}}ar is with the Coordinated Science Lab, University of Illinois
at Urbana-Champaign, Urbana, IL 61801, USA {\tt\small
basar1@illinois.edu}}
  }

\begin{document}

\maketitle
\thispagestyle{empty}
\pagestyle{empty}

\begin{abstract}


By enabling multiple agents to cooperatively solve a global
optimization problem  in the absence of a central coordinator,
decentralized stochastic optimization is gaining increasing
attention in areas as diverse as
 machine learning,  control, and sensor networks. Since the associated
 data usually contain sensitive information, such as user locations
 and personal identities, privacy protection has emerged as a crucial
 need in the implementation of decentralized stochastic
 optimization. In this paper, we propose a decentralized stochastic
 optimization algorithm that is able to guarantee provable
 convergence accuracy  even in the presence of aggressive
 quantization errors that are proportional to the amplitude of quantization inputs. The result applies to both convex and non-convex objective functions,
 and enables us to exploit aggressive quantization schemes to
obfuscate shared information, and hence enables privacy protection
without losing provable optimization accuracy. In fact, by using a {stochastic}
ternary quantization scheme, which quantizes any value to three
numerical levels, we achieve quantization-based rigorous
differential privacy in decentralized stochastic optimization, which
has not been reported before. In combination with the presented
quantization scheme, the proposed algorithm ensures, for the first
time, rigorous
 differential privacy in decentralized stochastic optimization without losing provable convergence
 accuracy. Simulation results for a
distributed estimation problem as well as numerical experiments for
decentralized   learning on a benchmark machine learning dataset
confirm the effectiveness of the proposed approach.

%
%
\end{abstract}

\section{Introduction}

{Initially introduced in the 1980s in the context of parallel and distributed computation \cite{tsitsiklis1984problems,bertsekas1989parallel}, decentralized optimization is finding increasing applications. For example, in sensor-network based acoustic-event localization, spatially distributed   sensors   multilaterate the position of a target event using individual sensors' range measurements such as time-of-arrival or signal-strength-profile measurements \cite{deligeorges2015mobile}. Because the range measurements acquired by individual sensors are noisy,  decentralized  optimization is commonly employed for the network to cooperatively   estimate the target position, particularly when the network is mobile or formed in an ad-hoc manner \cite{deligeorges2015mobile,zhang2017distributed}. Another example is the multi-robot rendezvous problem, where robots with different battery levels cooperatively determine a meeting time and place using decentralized optimization to minimize the total energy expenditure of the network \cite{cortes2006robust}. In wide-area monitoring and control of power systems, decentralized optimization enables multiple local control centers in a large power system network to cooperatively estimate and further damp inter-area electro-mechanical oscillations, which is vital for power system  stability \cite{nabavi2015distributed}. In large-scale machine learning, decentralized optimization algorithms are  becoming an important solution to parallelling both data and computation so as to handle the enormous growth in data and model sizes \cite{jiang2017collaborative}.}

 In decentralized
 optimization,   participating agents  interleave  on-device
 computation and peer-to-peer communications to cooperatively solve a network
 optimization problem. In recent years, a particular type of decentralized optimization, i.e., decentralized stochastic
 optimization, in which participating agents use  {\it noisy}
local gradients for optimization, is gaining increased traction due
to its superior performance in handling large or noisy data sets.
For example, in modern machine learning applications on massive
datasets,  such stochastic optimization methods are highly preferred
because they allow multiple devices to train a neural network model
collectively using local noisy gradients calculated from a small
batch of data points available to individual agents. Using a small
batch of data points yields a noisy estimation of the exact
gradient, but it is completely necessary because evaluating the
precise gradient using all available data can be extremely expensive
in computation or even practically infeasible. Furthermore, in the
era of Internet of things which connect  massive low-cost sensing
and communication devices, the data fed to optimization computations
are usually subject to measurement noises
\cite{xin2020decentralized}. As
 deterministic (batch) optimization approaches  typically falter
when dealing with   noisy data \cite{bottou2018optimization},
investigating  decentralized stochastic optimization algorithms
becomes a mandatory task.

Although centralized stochastic optimization algorithms can date
back to the 1950s \cite{bottou2018optimization}, results on
completely decentralized stochastic optimization in the absence of
any coordinator only started to gain attention in the past decade. So far,
plenty of decentralized stochastic optimization algorithms have been
reported, both for  convex objective functions (e.g.,
\cite{ram2010distributed,nedic2016stochastic,jakovetic2018convergence,sayin2017stochastic,pu2020distributed,rabbat2015multi,shamir2014distributed,sirb2018decentralized})
and non-convex objective functions (e.g.,
\cite{bianchi2012convergence,tatarenko2017non,lian2017can,singh2020sparq,koloskova2019decentralized,george2019distributed}).
In these decentralized stochastic optimization algorithms, because
participating agents only share gradients/model updates and do not
let raw data leave participants' machines, these algorithms were
believed to be able to protect the privacy of participating agents.
However, recent studies tell a completely different story: not only
can an adversary reversely infer the properties (e.g., membership
associations) of the data used in optimization
\cite{zhang2019admm,melis2019exploiting}, an adversary can even
precisely infer raw data used in optimization from shared gradients
(pixel-wise accurate for images and token-wise matching for texts)
\cite{zhu2019deep}. These information leakages pose a severe threat
to the privacy of participating agents in decentralized stochastic
optimization, as the data involved in optimization computation often
contain sensitive information such as medical records and financial
transactions.

Compared with  centralized optimization or distributed optimization
with a coordinator, privacy protection in  completely decentralized
  optimization is much more challenging due to the lack of
a trusted party. In fact, in decentralized stochastic optimization,
no participating agents are trustworthy as every participating agent
can use received messages to infer other participating agents'
sensitive information. Recently, results have been reported to
address the privacy issue in decentralized stochastic optimization.
One approach is to employ secure multi-party computation approaches
such as homomorphic encryption \cite{paillier1999public} or garbled
circuit \cite{yao1986generate}. However, {while allowing exact computations,} these approaches are very
heavy in computation/communication overhead, usually  incurring a
{\it runtime overhead of three to four orders of magnitude}
\cite{hynes2018efficient}. Furthermore, except our prior results
\cite{zhang2019admm,zhang2018enabling}, most existing homomorphic
encryption based privacy approaches employ  a  server (e.g., in
\cite{shoukry2016privacy,lu2018privacy,alexandru2020cloud}), which
does not exist in completely decentralized optimization. Hardware
based privacy approaches such as trusted hardware enclaves have also
been reported \cite{hynes2018efficient}. However, similar to
homomorphic encryption based approaches,  these approaches cannot be
directly used   to prevent multiple data providers from inferring
each others' data during decentralized stochastic optimization.
Another commonly used approach to enable privacy in decentralized
optimization is differential privacy, which adds uncorrelated noise
to shared gradients/model updates (e.g.,
\cite{bassily2014private,huang2015differentially,cortes2016differential,zhang2019recycled}).
However, these uncorrelated-noise based approaches are subject to a
fundamental trade-off between enabled privacy and optimization
accuracy \cite{dwork2014algorithmic}, i.e. a stronger privacy
protection requires a greater magnitude of uncorrelated noise, which
will unavoidably leads to a more intense reduction in optimization
accuracy.
 Recently,  results were reported to enable privacy by
exploiting the structural properties of decentralized optimization
\cite{li2020privacy,yan2012distributed,lou2017privacy}. For example,
the authors in
  \cite{yan2012distributed,lou2017privacy} proposed to add a {\it constant} uncertain parameter in projection or step
  sizes to enable privacy protection. The authors of \cite{gade2018privacy} proposed to judiciously construct spatially correlated ``structured" noise to cover
  gradient information without compromising optimization accuracy.
  However, the privacy protection enabled by these approaches is restricted: projection based privacy depends on
  the size of the projection set -- a large projection set nullifies
  privacy protection whereas a small projection set offers strong
  privacy
  protection but requires {\it a priori} knowledge of the optimal solution; ``structured"
  noise based approaches require  each agent to have a certain number
  of neighbors whose shared messages are inaccessible to  the adversary.  In
  fact, such a   constraint on information accessible to the adversary  is required in most existing
  accuracy-maintaining privacy solutions to decentralized optimization.
  For example, our studies in \cite{zhang2019admm} show that
   even  partially homomorphic encryption based privacy
  approaches  require  the adversary  not to have access to all messages shared by a target agent. One exception is our recent work \cite{wang2022decentralized,wang2022decentralized1}, which injects stochasticity in stepsizes and can  enable privacy protection even when adversaries have access to all shared messages in the network.

In this paper, we propose to leverage aggressive quantization
effects to enable strong privacy protection in decentralized
stochastic optimization without compromising optimization accuracy.
More specifically, we propose a decentralized stochastic
optimization algorithm that can ensure provable convergence accuracy
under aggressive quantization effects. This decentralized stochastic
optimization algorithm allows us to  quantize any shared value to
three numerical levels and hence obfuscate exchanged messages
without compromising optimization accuracy. In fact, we rigorously
prove that the quantization scheme can enable a strict
$(0,\delta)-$differential privacy for participating agents'
gradient information, which has not been reported in the literature.
The ability to use this aggressive quantization scheme also allows
us to significantly reduce communication overhead without losing
optimization accuracy since each real-valued message becomes
representable with two bits after quantization.

The main contributions of the paper are as follows: 1) We propose a
completely decentralized stochastic optimization algorithm that can
maintain provable optimization accuracy in the presence of
aggressive quantization errors that can be proportional to the norm
of input values. This is different from existing results that
require the quantization errors to be bounded
\cite{reisizadeh2019robust} or diminishing \cite{berahas2019nested}
with time. Furthermore, we obtain provable convergence for both
convex objective functions and non-convex objective functions, which
is different from \cite{reisizadeh2019exact} which only addresses
strongly convex objective functions; 2) We propose to use a {stochastic} ternary
quantization scheme to achieve rigorous $(0,\delta)-$differential
privacy, which has not been reported in the literature. Note that
$(0,\delta)-$differential privacy is stronger than the commonly used
$(\epsilon,\delta)-$differential privacy; 3) By integrating with
ternary quantization, our algorithm achieves rigorous
$(0,\delta)-$differential privacy under provable convergence
accuracy. To the best of  our knowledge, this is the first time both
rigorous $(0,\delta)-$differential privacy and provable convergence
accuracy are achieved simultaneously in decentralized stochastic
optimization; 4) The  ternary quantization scheme also enables us to
improve communication efficiency, which is crucial in scenarios
where the communication bandwidth is limited.


The paper is organized as follows:  Sec. II provides the problem
formulation. Sec. III presents the  decentralized stochastic
optimization algorithm. Sec. IV proves converge of all agents  to the same stationary point in the presence of aggressive quantization effects when
the objective functions are non-convex. Sec. V proves that the proposed
algorithm guarantees  convergence of all agents to the optimal solution in the presence of aggressive quantization effects
when the objective functions are  convex. Sec. VI proves that a
specific instantiation  of allowable quantization schemes can enable
rigorous $(0,\delta)-$differential privacy and, hence the proposed
algorithm can achieve rigorous $(0,\delta)-$differential privacy
with provable convergence accuracy. Sec. VII gives simulation
results as well as numerical experiments on a benchmark machine
learning dataset to confirm the obtained results. Finally Sec. VIII
concludes the paper.


{\bf Notation:} We use the symbol $\mathbb{R}$ to denote the set of
real numbers and $\mathbb{R}^d$ the Euclidean space of dimension
$d$. ${\bf 1}$ denotes a column vector  of appropriate dimension
 with all entries equal to 1. A vector is viewed as a column vector,
unless otherwise stated. For a vector $x$, $x_i$ denotes its $i$th
element. $A^T$ denotes the transpose of matrix $A$ and $x^Ty$
denotes the scalar product of two vectors $x$ and $y$. We use
$\langle\cdot\rangle$ to denote inner product  and $\|\cdot\|$ to
denote the standard Euclidean norm $\|x\|=\sqrt{x^T x}$. We use
$\|\cdot\|_1$ and $\|\cdot\|_{\infty}$ to denote the $\ell_1$ norm
$\|x\|_1=\sum_{i=1}^d|x_i|$ and the $\ell_\infty$ norm
$\|x\|_{\infty}=\max(|x_1|,|x_2|,\cdots,|x_d|)$, respectively. A
square matrix $A$ is said to be column-stochastic when its elements
in every column add up to one. A matrix $A$ is said to be
doubly-stochastic when both $A$ and $A^T$ are column-stochastic
matrices. We use $P(\mathcal{A})$ to denote the probability of an
event $\mathcal{A}$ and $\mathbb{E}\left[x\right]$ the
expected value  of a random variable $x$.

\section{Problem Formulation}


We  consider a network of $m$ agents solving the following
optimization problem cooperatively:
\begin{equation}\label{eq:optimization_formulation1}
\min\limits_{x\in\mathbb{R}^d} \frac{1}{m}\sum_{i=1}^m
f_i(x),\quad f_i(x)\triangleq\mathbb{E}_{\xi_i\sim
\mathcal{D}_i}\left[F_i(x,\xi_i)\right]
\end{equation}
where $x\in \mathbb{R}^d$ is the optimization variable common to all
agents but $F_i:\mathbb{R}^d\times \mathbb{R}\rightarrow\mathbb{R}$
is a local stochastic loss function private to agent $i$.
$\mathcal{D}_i$ is the local distribution of  data samples. In
practice, the distribution  $\mathcal{D}_i$ is usually unknown and
we only have access to $n_i$ realizations of it, { denoted by $\xi_{i,1}, \xi_{i,2}, \cdots, \xi_{i,n_i}$, where $\xi_{i,j}$ denotes the  $j$th random data sample of node $i$}. Thus $f_i(x)$ in
(\ref{eq:optimization_formulation1}) is usually determined by $
f_i(x)=\frac{1}{n_i}\sum_{j=1}^{n_i}F_i(x,\xi_{i,j}) $ which makes
(\ref{eq:optimization_formulation1}) the empirical risk minimization
problem. 

Because of the randomness in $F_i(x,\xi_i)$, the gradient that each
agent $i$ can obtain is subject to noises. We denote the gradient
that agent $i$ obtains at iteration $k$ for optimization as
$g_i^k(x,\xi_i)$, which will hereafter be abbreviated as $g_i^k$.
We make the following
standard assumption about $f_i(\cdot)$ and $g_i^k$:
\begin{Assumption 1}\label{Assumption:gradient}
\begin{enumerate}
 \item All   $f_i(\cdot)$ are Lipschitz continuous with Lipschitz
 gradients
 \[
 \|\nabla
f_i(x)-\nabla f_i(y)\|\leq L
  \|x-y\|,\: \forall x\in\mathbb{R}^d, y\in  \mathbb{R}^d,
  \] { and (\ref{eq:optimization_formulation1}) always has at least one optimal solution $x^{\ast}$, i.e., $\sum_{i=1}^m\nabla f_i(x^{\ast})=0$;}
  \item All $g_i^k$ satisfy
\[
\begin{aligned}
&\mathbb{E}_{\xi_i\backsim\mathcal{D}_i}\left[g_i^k\right]=\nabla
f_i(x_i^k),\:\forall i &\\
&\mathbb{E}_{\xi_i\backsim\mathcal{D}_i}\left[\|g_i^k-\nabla
  f_i(x)\|^2\right]\leq \sigma^2,\: \forall i, x&
\end{aligned}
\]
\end{enumerate}
\end{Assumption 1}

In order for the network of $m$ agents to cooperatively solve
(\ref{eq:optimization_formulation1})  in a decentralized manner, we
assume that the $m$ agents interact on an undirected graph. The
interaction can be described by a weight matrix $W$. More
specifically, if agent $i$ and agent $j$ can communicate and
interact with each other, then the $(i,j)$th entry of $W$, i.e.,
$w_{ij}$, is positive. Otherwise, $w_{ij}$ is zero. 
 The neighbor set
$\mathcal{N}_i$ of agent $i$ is defined as the set of agents
satisfying $\{j|w_{ij}>0\}$.  We define a diagonal matrix $D$ with
the $i$th diagonal entry determined  as
$d_{ii}=\sum_{j\in{\mathbb{N}}_i}w_{ij}$. So the matrix $D-W$ will
be the commonly referred graph Laplacian matrix.   To ensure that
the network can cooperatively solve
(\ref{eq:optimization_formulation1}), we make the following standard
assumption about the interaction:

\begin{Assumption 2}\label{assumption:coupling}
   The interaction topology forms an undirected connected network,
   i.e.,  the second smallest eigenvalue $\rho$ of the graph Laplacian matrix $L_w\triangleq D-W$
 is positive.
\end{Assumption 2}

In decentralized stochastic optimization, gradients are directly
computed from raw data and hence  embed   sensitive information. For
example, in decentralized-optimization based
 localization, disclosing the gradient  of an agent
amounts to disclosing its  position
\cite{zhang2019admm,huang2015differentially}. In machine learning
applications, gradients are directly calculated  from and embed
information of sensitive training data \cite{zhu2019deep}.
Therefore, in this paper, we define privacy as preventing  agents'
gradients from being inferable by adversaries.

  We consider two potential adversaries in decentralized stochastic optimization, which are the
two most commonly used models of attacks in privacy research
\cite{Goldreich_2} and have been widely used in the control community \cite{wang2019privacy,ruan2019secure}:
\begin{itemize}
\item \emph{Honest-but-curious attacks}  are attacks in which a participating
agent or multiple participating agents (colluding or not)
 follows all protocol steps correctly but is curious and
collects all received intermediate data in an attempt to learn the
sensitive information about other participating agents.

\item \emph{Eavesdropping attacks}  are attacks in which an external eavesdropper wiretaps all communication channels to
intercept exchanged messages so as to learn sensitive information
about sending agents.
\end{itemize}

{ An honest-but-curious adversary (e.g., agent $i$) has access to the internal state $x_i^k$, which is unavailable
to external eavesdroppers. However, an eavesdropper has access to all shared information in the network,
whereas an honest-but-curious agent can only access shared information that is destined to it.}

In this paper, we propose to leverage quantization effects to enable
differential privacy in decentralized stochastic optimization. We
adopt the definition of ($\epsilon,\delta$)-differential privacy
following standard conventions \cite{dwork2014algorithmic}:
\begin{Definition 1}
For a randomized function $h(x)$, we say that it is $(\epsilon,
\delta)$-differentially private if for all { subsets $S$ of the image set of the function $h(x)$} 
and for all $x,y$ with $\|x-y\|_1\leq 1$, we always have
\[
P(h(x)\in S)\leq e^\epsilon P(h(y)\in S)+\delta.
\]
\end{Definition 1}

{ Definition 1 says that for two inputs $x$ and $y$ with $\ell_1$-norm difference no more than $1$, a mechanism $h(\cdot)$ achieves $(\epsilon, \delta)-$differential privacy if it can ensure that the outputs of the two inputs are different in   probabilities by at most   $\epsilon$ and $\delta$ specified on the right hand side of the above inequality. Clearly, a smaller $\epsilon\geq 0$ or $\delta\geq 0$ means   better differential-privacy protection. In Sec. VI we will prove that a specific quantization mechanism can enable $(0, \delta)-$differential privacy protection for  exchanged information. {Note that under a fixed value of $\delta$, $(0, \delta)-$differential privacy  is  stronger than $(\epsilon, \delta)-$differential privacy for any $\epsilon> 0$.}

\begin{Remark 1}
   In the original definition of differential privacy in \cite{dwork2014algorithmic,meiser2018approximate}, because the  input space is discrete, i.e., $x$ and $y$ are strings, the distance between $x$ and $y$ is measured by the number of positions at which the corresponding symbols are different  (Hamming distance).   In our case, since the input space is continuous, we use $\ell_1$ norm to measure the distance between two real vectors $x$ and $y$. In fact, any $\ell_p$ norm defined by $\|x\|_p=\left(|x_1|^p+|x_2|^p+\cdots+|x_m|^p\right)^{1/p}$ with $p\geq 1$ can be used in the definition.
\end{Remark 1}
}

\section{Quantization-enabled privacy-preserving decentralized optimization algorithm}

\begin{figure}
\centering
\begin{minipage}{.49\textwidth}
\small
\par\noindent\rule{\textwidth}{0.5pt} \noindent\textbf{Algorithm 1:
Quantization-enabled Privacy-preserving Decentralized Stochastic
Optimization}

\vspace{-0.2cm}\noindent\rule{0.99\textwidth}{0.5pt}
\begin{enumerate}
\item Public parameters: $W$, $\epsilon^k$, $\lambda^k$ $x_i^0=0$ for all $i$,
the total number of iterations $t$
\item For the $i$th agent, at iteration $k$
\begin{enumerate}
\item Determine local gradient $g_i^k$;
\item Determine quantized state $\mathcal{Q}(x_i^k)$ and send
it to all agents $j\in\mathbb{N}_i$;
\item After receiving $\mathcal{Q}(x_{j}^k)$ from all
$j\in\mathbb{N}_i$, update state as
\[
\begin{aligned}
x_i^{k+1}&= x_i^k+\epsilon^k\sum_{j\in\mathbb{N}_i} w_{ij}
(\mathcal{Q}(x_j^k)-\mathcal{Q}(x_i^k) )-\epsilon^k\lambda^k g_i^k
\end{aligned}
\]
\end{enumerate}
    \item end
\end{enumerate}
\vspace{-0.2cm}\rule{0.99\textwidth}{0.5pt}
\end{minipage}
\end{figure}

Before presenting our quantization-enabled privacy-preserving
approach for decentralized stochastic optimization, we first discuss
why conventional decentralized stochastic optimization algorithms
leak gradient information of participating agents.

{By assigning a copy $x_i$ of the decision variable $x$ to each agent $i$, and then imposing the requirement  $x_i = x$ for all $1\leq i\leq m$, we can rewrite the
optimization problem (\ref{eq:optimization_formulation1}) in the following form \cite{nedic2020distributed}:}
\begin{equation}
\min\limits_{x\in\mathbb{R}^{md}} f(x)= \frac{1}{m}\sum_{i=1}^m
f_i(x_i) \quad {\rm s.t.} \quad x_1=x_2=\cdots=x_m
\end{equation}
where $x=[x_1^T,x_2^T,\cdots,x_m^T]^T$.  Conventional decentralized
optimization algorithms usually take the following form
\cite{jiang2017collaborative,lian2017can}:
\[
x_i^{k+1}=x_i^k+ \sum_{j\in\mathbb{N}_i}  w_{ij}(x_j^k-x_i^k)-\eta
g_i^k
\]
where $x_i^k$ denotes the optimization variable maintained by agent
$i$ at iteration $k$, and $\eta$ denotes the optimization stepsize, { which should be no greater than $\frac{1}{L}$ to ensure stability \cite{lian2017can}}.
Because $w_{ij}$ has to be publicly known to establish conditions in
Assumption \ref{assumption:coupling} in a decentralized manner
\cite{gharesifard2012distributed} and
 agent $i$  shares $x_i^k$ with all its neighbors,  an
adversary can calculate the gradient $g_i^k$ of any agent based on
publicly known $W$ and $\eta$ if it   has access to all information
shared in the network.

Motivated by this observation, we propose the following
decentralized optimization algorithm which leverages quantization to
enable privacy protection:
\begin{equation}\label{eq:proposed_algorithm}
\begin{aligned}
x_i^{k+1}&= x_i^k+\epsilon^k\sum_{j\in\mathbb{N}_i}
w_{ij}\left(\mathcal{Q}(x_j^k)-\mathcal{Q}(x_i^k)\right)-\epsilon^k\lambda^k
g_i^k
\end{aligned}
\end{equation}
where $\lambda^k$ and $\epsilon^k$ are publicly-known design
parameters   crucial for ensuring provable convergence accuracy
under aggressive quantization effects, and their design will be
elaborated on later. Note that, although agent $i$ has access to
 $x_i^k$, we still use a quantized version of $x_i^k$ in the
comparison term $\mathcal{Q}(x_j^k)-\mathcal{Q}(x_i^k)$ in
(\ref{eq:proposed_algorithm}). This is intuitive as when $x_i^k$ and
$x_j^k$ are the same, we do not want the quantization operation to
introduce an extra non-zero input to the optimization process. In
fact, as shown in later derivations, this strategy will also
simplify the evolution of the average optimization variable across
all agents.

In our proposed algorithm (\ref{eq:proposed_algorithm}), at
iteration $k$,  every  agent $i$ only shares quantized state $x_i^k$
(see details in Algorithm 1). Therefore, even if an adversary has
access to the quantized state of an agent $i$ as well as all
information received by agent $i$ (which are also quantized), the
adversary still cannot use the dynamics
(\ref{eq:proposed_algorithm}) to precisely infer the gradient of
agent $i$ due to  quantization induced errors. In fact, as will be
proved later, the proposed algorithm can
 have provable convergence even in the presence of aggressive
quantization schemes with large quantization errors, which will
enable us to achieve strict  $(0,\delta)$-differential privacy
  protection for all participating agents. {More specifically,
we consider stochastic quantization schemes satisfying the following
Assumption}:

\begin{Assumption 1}\label{Assmption:quantization}
The   quantizer $\mathcal{Q}(\cdot)$ is unbiased and its variance is
proportionally bounded by the input's norm, i.e., $
\mathbb{E}\left[\mathcal{Q}(x)|x\right]=x$ and $
\mathbb{E}\left[\|\mathcal{Q}(x)-x\|^2|x\right]\leq \beta\|x\|^2$
hold for some constant $\beta$ and any $x$. And the quantization on
different agents are independent of each other.
\end{Assumption 1}

\begin{Remark 1}
Note that the quantization schemes considered in Assumption
\ref{Assmption:quantization} are quite general and include the
commonly used error-bounded quantization schemes (in, e.g.,
\cite{rabbat2005quantized,yi2014quantized,reisizadeh2019robust}) and
error-diminishing quantization schemes (in, e.g.,
\cite{pu2016quantization,berahas2019nested}) as special cases.
\end{Remark 1}

\begin{Remark 1}
Note that when the quantization scheme is designed such that it only
outputs the sign of the quantization input (which still satisfies
the conditions in Assumption \ref{Assmption:quantization}), the
inter-agent coupling in the proposed algorithm looks similar to the
interaction in existing decentralized optimization algorithms that
use only  the sign of relative states (see,
\cite{zhang2019distributed,cao2021decentralized}). However, there is
a crucial difference between the two in that the  quantization
scheme here can be implemented by every participating agent without
knowing anything about its neighbors' states, whereas the
relative-state sign based interaction (which arises in other
contexts) requires an agent to know (some) information about its
neighbors' states.
\end{Remark 1}

Augmenting the decision variables of all agents as
$x^k=[(x_1^k)^T,(x_2^k)^T,\cdots,(x_m^k)^T]^T$,
  we can write the overall network dynamics of
the proposed decentralized optimization algorithm as follows
\begin{equation}\label{eq:proposed_algorithm_vector_form}
x^{k+1}=(A^k\otimes I_d) x^k-\epsilon^k \lambda^k g^k-\epsilon^k
(L_w\otimes I_d) V^k
\end{equation}
where $L_w$ is the Laplacian matrix defined in Assumption
\ref{assumption:coupling},
\[
A^k=(I-\epsilon^kL) \in\mathbb{R}^{m \times  m},
\]
\[
g^k=\left[(g_1^k)^T,\,(g_2^k)^T,\,\cdots,(g_m^k)^T\right]^T\in
\mathbb{R}^{md\times1},
\]
\[
V^k=\left[(v_1^k)^T,\,(v_2^k)^T,\,\cdots,(v_m^k)^T\right]^T\in
\mathbb{R}^{md\times1},
\]
\[
v_i^k=  \mathcal{Q}(x_i^k)-x_i^k\in \mathbb{R}^{d\times1}
\] Here $\otimes$
denotes Kronecker product and $I_d$ denotes identity matrix of
dimension $d$.

It can be obtained that the evolution of the average optimization
variable $\bar{x}^k=\frac{\sum_{i=1}^mx_i^k}{m}$ follows
\begin{equation}\label{eq:bar_x_evolution}
\begin{aligned}
\bar{x}^{k+1}&=\bar{x}^k+\frac{\epsilon^k}{m}\sum_{i=1}^{m} \sum_{j\in\mathbb{N}_i}
w_{ij}\left(\mathcal{Q}(x_j^k)-\mathcal{Q}(x_i^k)\right)\\
&\qquad\qquad\qquad\qquad\qquad\qquad-\epsilon^k\lambda^k\frac{\sum_{i=1}^mg_i^k}{m}\\
&=\bar{x}_i^k-\epsilon^k\lambda^k\frac{\sum_{i=1}^mg_i^k}{m}
\end{aligned}
\end{equation}
which is independent of the quantization error. { Note that in the second equality, we used the fact that the network is undirected, i.e.,   $w_{ij}=w_{ji}$ from Assumption \ref{assumption:coupling}, which leads to the annihilation of all coupling terms due to $w_{ij}\left(\mathcal{Q}(x_j^k)-\mathcal{Q}(x_i^k)\right)+w_{ji}\left(\mathcal{Q}(x_i^k)-\mathcal{Q}(x_j^k)\right)=0$.} This shows the
benefit for agent $i$ to use its quantized state $x_i^k$ in  the
comparison term $\mathcal{Q}(x_j^k)-\mathcal{Q}(x_i^k)$ on the right
hand side of (\ref{eq:proposed_algorithm}).

{
\begin{Remark 1}
  From the above argument, it can be seen that  agents being able to update in a synchronized manner is key to guaranteeing the average optimization variable $\bar{x}^k$ to be immune to aggressive quantization errors.
\end{Remark 1}
}

In  the following two sections, we will show that   the proposed
decentralized stochastic optimization algorithm still has provable
convergence accuracy under aggressive quantization effects. More
specifically, in Sec. \ref{section:non_convex}, we will show that in the
non-convex case, the algorithm guarantees provable convergence of all
agents to the  same stationary
point; in Sec. \ref{section:convex}, we
will show that in the  convex case, the algorithm
guarantees the convergence of all agents to the optimal solution.


\section{Convergence  Analysis In the Non-convex
Case}\label{section:non_convex}

In this section, we show that the
proposed algorithm will ensure convergence of all agents   to the
same stationary point  when the objective functions are non-convex,
even under aggressive quantization effects.

To this end, we first show that when $\epsilon^k$ and $\lambda^k$ are
chosen appropriately, $\|g_i^k\|$ and
$\mathbb{E}\left[\|x^k\|^2\right]$ will always be bounded, which
allows us to quantify the effects of quantization on the
optimization process ({note that here the expectation is taken with respect to the randomness in stochastic gradients  and quantization up until iteration $k-1$}).  It is worth noting that as the results are
obtained irrespective of the convexity of objective  functions, they
are applicable to the derivations in the  convex case in the next
section, too.

\begin{Lemma 1}\label{Lemma:bounded_gradient}
Under Assumption \ref{Assumption:gradient}, the gradient $\|g_i^k\|$
is always bounded by some constant $G$.
\end{Lemma 1}
\begin{proof}
Under the conditions in Assumption \ref{Assumption:gradient}, the
result can be easily obtained from \cite{george2019distributed} or
Lemma 3.3 in \cite{khalil2002nonlinear}.
\end{proof}

\begin{Lemma 2}\label{lemma:lemma_boundedx^2}
Under Assumption \ref{Assumption:gradient}, Assumption
\ref{assumption:coupling}, and Assumption
\ref{Assmption:quantization}, $\mathbb{E}\left[\|x^k\|^2\right]$
will always be bounded if the positive sequences $\epsilon^k$ and
$\lambda^k$ satisfy $\sum_{k=1}^{\infty}(\epsilon^k)^2<\infty$ and
$\sum_{k=1}^{\infty}\epsilon^k(\lambda^k)^2<\infty$, { where the expectation is taken with respect to the randomness in stochastic gradients  and quantization up until iteration $k-1$.}
\end{Lemma 2}

\begin{proof}
The proof is given in Appendix B.
\end{proof}

Using Lemma \ref{lemma:lemma_boundedx^2}, we can further obtain that
the optimization variables $x_i^k$ of different agents will converge
to the average optimization variable  across all agents $\bar{x}^k$:
\begin{Lemma 3}\label{Lemma:convergence_of_x_i_to_bar_x}
Under the conditions in Lemma \ref{lemma:lemma_boundedx^2}, the
proposed algorithm guarantees
\[
\lim_{k\rightarrow \infty}
\mathbb{E}\left[\|x^{k+1}-\hat{\bar{x}}^{k+1}\|^2\right]=0
\]
where
$\hat{\bar{x}}^k\triangleq
 {\bf 1}_m\otimes\bar{x}^k$ with  ${\bf 1}_m$ denoting the $m$ dimensional column vector
of  $1$s. More specifically, represent the decaying rate of $
\lambda^k$   and $\epsilon^k$ as $0<\delta_1<1$ and $0<\delta_2<1$,
respectively, i.e., there exist some positive $a_1$, $a_2$, and $a_3$ such that
$\lambda^k \leq \frac{a_1}{(a_3k+1)^{\delta_1}}$    and $\epsilon^k
\leq \frac{a_2}{(a_3k+1)^{\delta_2}}$ hold,  then we have
\[
\lim_{k\rightarrow \infty}(1+k)^\delta
\mathbb{E}\left[\|x^{k+1}-\hat{\bar{x}}^{k+1}\|^2\right]=0
\]
for
any $0\leq \delta <\min\{2\delta_1, \delta_2\}$.
\end{Lemma 3}
\begin{proof}
The proof is given in Appendix C.
\end{proof}

Based on these
results, we can prove the following results on the convergence of
all agents to the same stationary point where the gradients are
zero:

\begin{Theorem 1}\label{theorem_convergence_nonconvex}
Under  Assumptions \ref{Assumption:gradient},
\ref{assumption:coupling}, and \ref{Assmption:quantization}, when the sequences
$\epsilon^k$ and $\lambda^k$ are selected such that the sequence
$\epsilon^k\lambda^k$ is not summable, but  $(\epsilon^k)^2$ and $
\epsilon^k(\lambda^k)^2$   are summable, i.e.,
\begin{equation}\label{eq:conditions}
\sum_{k=1}^{\infty} \epsilon^k \lambda^k=+\infty,\:\:
\sum_{k=1}^{\infty}(\epsilon^k)^2<\infty,\:\: \sum_{k=1}^{\infty}
\epsilon^k(\lambda^k)^2<\infty
\end{equation}
 then
    the proposed  algorithm will guarantee the following results:
\begin{equation}\label{eq:gradient_summable1}
\begin{aligned}
&\lim_{t\rightarrow\infty}\frac{\sum_{k=0}^t
\epsilon^k\lambda^k\mathbb{E}\left[\left\|\nabla
f(\bar{x}^k)\right\|^2\right]}{\sum_{k=0}^t \epsilon^k \lambda^k}=0,\\
&\lim_{t\rightarrow\infty}\frac{\sum_{k=0}^t \epsilon^k\lambda^k
\mathbb{E}\left[\left\|\frac{\sum_{i=1}^m\nabla f_i(x_i^k)
}{m}\right\|^2\right]}{\sum_{k=0}^t \epsilon^k\lambda^k}=0
\end{aligned}
\end{equation}
{ where the expectation is taken with respect to the randomness in stochastic gradients  and quantization up until iteration $k-1$.}
\end{Theorem 1}

\begin{proof}
From the Lipschitz gradient condition  in Assumption
\ref{Assumption:gradient}, we have
\[
f(y)\leq f(x)+ \langle \nabla f(x),y-x\rangle +\frac{L\|y-x\|^2}{2}
\]
for any $x\in\mathbb{R}^d$ and  $y\in\mathbb{R}^d$. By plugging
$y=\bar{x}^{k+1}$ and $x=\bar{x}^{k}$ into the above inequality, we
can have the following relationship based on
(\ref{eq:bar_x_evolution}):
\begin{equation}
\begin{aligned}
f(\bar{x}^{k+1})&\leq f(\bar{x}^k)+\left\langle \nabla f(\bar{x}^k),
-\epsilon^k\lambda^k\frac{\sum_{i=1}^mg_i^k}{m}\right\rangle
\\
&\qquad+\frac{L}{2}\left\|-\epsilon^k\lambda^k\frac{\sum_{i=1}^mg_i^k}{m}\right\|^2
\end{aligned}
\end{equation}
Taking expectation on both sides, we  can obtain

\begin{equation}\label{eq:f_X_k+1_secondstep}
\begin{aligned}
&\mathbb{E}\left[f(\bar{x}^{k+1})\right]\\
&\leq
\mathbb{E}\left[f(\bar{x}^k)\right]+\mathbb{E}\left[\left\langle
\nabla
f(\bar{x}^k),-\epsilon^k\lambda^k\frac{\sum_{i=1}^mg_i^k}{m}\right\rangle
\right]\\
&\qquad+
\frac{L}{2}\mathbb{E}\left[\left\|-\epsilon^k\lambda^k\frac{\sum_{i=1}^mg_i^k}{m}\right\|^2\right]\\
&=
\mathbb{E}\left[f(\bar{x}^k)\right]-\epsilon^k\lambda^k\mathbb{E}\left[\left\langle
\nabla f(\bar{x}^k),\frac{\sum_{i=1}^m g_i^k}{m}\right\rangle
\right]\\
&\qquad+
\frac{L(\epsilon^k\lambda^k)^2}{2m^2}\mathbb{E}\left[\left\|-
\sum_{i=1}^m g_i^k  \right\|^2\right]
\end{aligned}
\end{equation}

Using the equality $2\langle X,Y\rangle=\|X\|^2+\|Y\|^2-\|X-Y\|^2$,
we arrive at  the following relationship for the second term on the
right hand side of (\ref{eq:f_X_k+1_secondstep}):
\begin{equation}\label{eq:inner_product2}
\begin{aligned}
&\mathbb{E}\left[\left\langle  \nabla  f(\bar{x}^k), \frac{
\sum_{i=1}^m  g_i^k}{m} \right\rangle\right]\\
& =\mathbb{E}\left[\left\langle
\nabla f(\bar{x}^k),  \frac{\sum_{i=1}^m \nabla f_i(x_i^k)}{m}\right\rangle\right]\\
&= \frac{1}{2}\mathbb{E}\left[\left\| \nabla
f(\bar{x}^k)\right\|^2\right] +
 \frac{1}{2} \mathbb{E}\left[\left\|\frac{\sum_{i=1}^m\nabla
f_i(x_i^k) }{m}\right\|^2\right]\\
&\qquad - \frac{1}{2}\mathbb{E}\left[\left\|\nabla f(\bar{x}^k)
-\frac{ \sum_{i=1}^m \nabla
f_i(x_i^k)}{m}\right\|^2 \right] \\
&\geq   \frac{1}{2}\mathbb{E}\left[\left\| \nabla
f(\bar{x}^k)\right\|^2\right] +
  \frac{1}{2}\mathbb{E}\left[\left\|\frac{\sum_{i=1}^m\nabla
f_i(x_i^k) }{m}\right\|^2\right]\\
&\qquad -\frac{ L^2}{2m}\sum_{i=1}^m
\mathbb{E}\left[\left\|\bar{x}^k - x_i^k\right\|^2 \right]
\end{aligned}
\end{equation}
where we used the   Lipschitz gradient assumption in Assumption
\ref{Assumption:gradient} and the relationship
$\|y_1+y_2+\cdots+y_m\|^2\leq m\sum_{i=1}^m\|y_i\|^2$ in the
inequality.

For the third term on the right hand side of
(\ref{eq:f_X_k+1_secondstep}), we can bound it using the result that
$g_i^k$ is bounded by $G$ obtained in  Lemma
\ref{Lemma:bounded_gradient}:

\begin{equation}\label{eq:third_term}
\begin{aligned}
\frac{L(\epsilon^k\lambda^k)^2}{2m^2}\mathbb{E}\left[\left\|-
\sum_{i=1}^m g_i^k  \right\|^2\right] &\leq
\frac{L(\epsilon^k\lambda^k)^2}{2m} \mathbb{E}\left[\sum_{i=1}^m \left\|g_i^k\right\|^2\right]  \\
&\leq \frac{LG^2(\epsilon^k\lambda^k)^2}{2 }
\end{aligned}
\end{equation}

Plugging (\ref{eq:inner_product2}) and (\ref{eq:third_term}) into
(\ref{eq:f_X_k+1_secondstep}) leads to

\begin{equation}\label{eq:f_X_k+1_thirdstep}
\begin{aligned}
\mathbb{E}\left[f(\bar{x}^{k+1})\right]&\leq
\mathbb{E}\left[f(\bar{x}^k)\right]-\frac{1}{2}\epsilon^k\lambda^k
\mathbb{E}\left[\left\| \nabla f(\bar{x}^k)\right\|^2\right]\\
&\quad -\frac{1}{2}\epsilon^k\lambda^k
  \mathbb{E}\left[\left\|\frac{\sum_{i=1}^m\nabla
f_i(x_i^k) }{m}\right\|^2\right]\\
&\quad+\epsilon^k\lambda^k\frac{ L^2}{2m}\sum_{i=1}^m
\mathbb{E}\left[\left\|\bar{x}^k  - x_i^k\right\|^2
\right]\\
&\quad +\frac{LG^2(\epsilon^k\lambda^k)^2}{2 }
\end{aligned}
\end{equation}
or
\begin{equation}\label{eq:theorem2}
\begin{aligned}
&\epsilon^k\lambda^k \mathbb{E}\left[\left\| \nabla
f(\bar{x}^k)\right\|^2\right] +\epsilon^k\lambda^k
  \mathbb{E}\left[\left\|\frac{\sum_{i=1}^m\nabla
f_i(x_i^k) }{m}\right\|^2\right]\\
&\leq2\left(\mathbb{E}\left[f(\bar{x}^k)\right]-\mathbb{E}\left[f(\bar{x}^{k+1})\right]\right)\\
&\quad+\epsilon^k\lambda^k\frac{ L^2}{m}\sum_{i=1}^m
\mathbb{E}\left[\left\|\bar{x}^k  - x_i^k\right\|^2 \right]+
LG^2(\epsilon^k\lambda^k)^2
\end{aligned}
\end{equation}

Iterating the above inequality from $k=0$ to $k=t$ yields
\begin{equation}
\begin{aligned}
\sum_{k=0}^{t} &\left(\epsilon^k\lambda^k \mathbb{E}\left[\left\|
\nabla f(\bar{x}^k)\right\|^2\right]  +\epsilon^k\lambda^k
  \mathbb{E}\left[\left\|\frac{\sum_{i=1}^m\nabla
f_i(x_i^k) }{m}\right\|^2\right]\right) \\
&\leq
 2\left(\mathbb{E}\left[f(\bar{x}^{0})\right]-\mathbb{E}\left[f(\bar{x}^{t+1})\right]\right)\\
&  \quad +\sum_{k=0}^{t}\epsilon^k\lambda^k\frac{L^2}{m}\sum_{i=1}^m
\mathbb{E}\left[\left\|\bar{x}^k  - x_i^k\right\|^2 \right]\\
&\quad  +\sum_{k=0}^{t} LG^2(\epsilon^k\lambda^k)^2
\end{aligned}
\end{equation}
i.e.,
\begin{equation}\label{eq:final}
\begin{aligned}
&\frac{\sum_{k=0}^{t}\epsilon^k \lambda^k\mathbb{E}\left[
\left\|\nabla
f(\bar{x}^k)\right\|^2\right]}{\sum_{k=0}^{t}\epsilon^k \lambda^k}\\
&\qquad  +\frac{\sum_{k=0}^{t}\epsilon^k  \lambda^k\mathbb{E}\left[
\left\|\frac{\sum_{i=1}^m\nabla f_i(x_i^k)
}{m}\right\|^2\right]}{\sum_{k=0}^{t}\epsilon^k\lambda^k}
\\
& \leq
 \frac{2\left(\mathbb{E}\left[f(\bar{x}^{0})\right]-\mathbb{E}\left[f(\bar{x}^{t+1})\right]\right)}{\sum_{k=0}^{t}\epsilon^k \lambda^k}\\
&  \qquad
+\frac{\sum_{k=0}^{t}\epsilon^k\lambda^k\frac{L^2}{m}\sum_{i=1}^m
\mathbb{E}\left[\left\|\bar{x}^k  - x_i^k\right\|^2
\right]}{\sum_{k=0}^{t}\epsilon^k \lambda^k}\\
&\qquad  +\frac{\sum_{k=0}^{t}
LG^2(\epsilon^k\lambda^k)^2}{\sum_{k=0}^{t}\epsilon^k\lambda^k}
\end{aligned}
\end{equation}
It can be verified that when
 $\epsilon^k$ and $\lambda^k$ are selected in such a way
that the conditions in (\ref{eq:conditions}) are satisfied, then the
conditions in Lemma \ref{Lemma:convergence_of_x_i_to_bar_x} will
also be satisfied, which means that
$\mathbb{E}\left[\left\|\bar{x}^k  - x_i^k\right\|^2 \right]$ will
be in the same order as $(\lambda^k)^2$
 or $\epsilon^k$. This means that
$\epsilon^k\lambda^k\frac{ L^2}{m}\sum_{i=1}^m
\mathbb{E}\left[\left\|\bar{x}^k  - x_i^k\right\|^2 \right]$ will be
in the  same order as $\epsilon^k( \lambda^k)^3$
  or
$(\epsilon^k)^2 \lambda^k$, both of which are summable according to
the conditions in (\ref{eq:conditions}). Therefore, the second term
on the right hand side of (\ref{eq:final}) will converge to zero.
Similarly, we can prove that all other terms on the right hand side
of (\ref{eq:final}) will converge to zero under the conditions in
(\ref{eq:conditions}), which completes the proof.
\end{proof}

{
\begin{Remark 1}
 Using the   Stolz-Ces\`{a}ro theorem, one can obtain from (\ref{eq:gradient_summable1}) that the limit inferiors of $\mathbb{E}\left[\|\nabla f(\bar{x}^t)\|^2\right]$ and $\mathbb{E}\left[\|\nabla f_i(x_i^t)\|^2\right]$ are zero as $t$ tends to infinity, i.e., $\underline{\lim}_{t\rightarrow\infty}\mathbb{E}\left[\|\nabla f(\bar{x}^t)\|^2\right]=0$ and $\underline{\lim}_{t\rightarrow\infty}\mathbb{E}\left[\|\nabla f_i(x_i^t)\|^2\right]=0$.
\end{Remark 1}}

In fact, if we can specify the convergence rate of $\epsilon^k$ and
$\lambda^k$, we can further obtain the convergence rate of the algorithm:

\begin{Corollary 1}\label{Theorem:convergence_rate}
 If the sequences $\epsilon^k$ and $\lambda^k$ are selected in the form of $\lambda^k
=\frac{a_{1}}{(a_3k+1)^{\delta_1}}$ and $\epsilon^k =
\frac{a_2}{(a_3k+1)^{\delta_2}}$ with $a_1$, $a_2$, and $a_3$ denoting some
positive constants and  positive exponents $\delta_1$ and $\delta_2$
satisfying $\delta_1+\delta_2\leq1$, $\delta_2>0.5$, and
$2\delta_1+\delta_2>1$,  then all conditions in
(\ref{eq:conditions}) are satisfied and the proposed algorithm will
guarantee (\ref{eq:gradient_summable1}) under Assumptions
\ref{Assumption:gradient}, \ref{assumption:coupling}, and
\ref{Assmption:quantization}. More specifically, the
convergence rate of gradients satisfies
\begin{equation}\label{eq:rate_in_nonconvex1}
\begin{aligned}
&\frac{\sum_{k=0}^{t}\epsilon^k \lambda^k\mathbb{E}\left[
\left\|\nabla
f(\bar{x}^k)\right\|^2\right]}{\sum_{k=0}^{t}\epsilon^k \lambda^k}
\\
&\qquad+\frac{\sum_{k=0}^{t}\epsilon^k  \lambda^k\mathbb{E}\left[
\left\|\frac{\sum_{i=1}^m\nabla f_i(x_i^k)
}{m}\right\|^2\right]}{\sum_{k=0}^{t}\epsilon^k\lambda^k}\\
&\qquad+
\frac{2\left(\mathbb{E}\left[f(\bar{x}^{t+1})\right]-\mathbb{E}\left[f(\bar{x}^{0})\right]\right)}{\sum_{k=0}^{t}\epsilon^k
\lambda^k}\\
& =\mathcal{O}\left(\frac{1}{(t+1)^{\delta}}\right)
\end{aligned}
\end{equation}
where $\delta=\min\{2\delta_1,\delta_2\}$ and { the expectation is taken with respect to the randomness in stochastic gradients  and quantization up until iteration $k-1$}.
\end{Corollary 1}
\begin{proof}
The proof  follows from the line of derivation in the proof of
Theorem \ref{theorem_convergence_nonconvex}. More specifically,
under the conditions of Theorem \ref{Theorem:convergence_rate}, the
conditions of Lemma \ref{Lemma:convergence_of_x_i_to_bar_x} will be
satisfied and we have the second term on the right hand side of
(\ref{eq:final}) converging to zero with a   rate of no less than
$\mathcal{O}\left(\frac{1}{(t+1)^{\delta}}\right)$ with $\delta
=\min\{2\delta_1, \delta_2\}$. Further note that the last term on
the right hand side of (\ref{eq:final}) converges to zero with a
rate $\mathcal{O}\left(\frac{1}{(t+1)^{\delta}}\right)$ with $\delta
=\delta_1+\delta_2$. Therefore,
 we have that the left hand side of (\ref{eq:rate_in_nonconvex1})
  will decay with a rate
$\delta=\min\{2\delta_1,\delta_2\}$ as defined in the statement.
\end{proof}

\section{Convergence  Analysis In the Convex
Case}\label{section:convex}

In this section, we consider the case where the objective functions are convex:
\begin{Assumption 2}\label{Assumtion:convex}
The objective functions $f_i(\cdot)$ are convex.
\end{Assumption 2}

As the derivations of the results in Lemma \ref{lemma:lemma_boundedx^2} and Lemma
\ref{Lemma:convergence_of_x_i_to_bar_x} are independent of the convexity of $f_i(\cdot)$, we still have the same results in the convex case. Therefore, in the convex case we can still have the same results obtained in Theorem \ref{theorem_convergence_nonconvex}. Moreover,  we can   prove that the convexity assumption in Assumption \ref{Assumtion:convex} also enables us to characterize convergence in function value to the optimal solution:
{
\begin{Theorem 1}\label{Theorem:convergence_convex}
  Under Assumptions \ref{Assumption:gradient}-\ref{Assumtion:convex},
when the positive sequences $\epsilon^k$ and $\lambda^k$ are
selected
such that the sequence
$\epsilon^k\lambda^k$ is not summable, but  $(\epsilon^k)^2$ and $
\epsilon^k(\lambda^k)^2$   are summable, i.e.,
\begin{equation}\label{eq:conditions_convex}
\sum_{k=1}^{\infty} \epsilon^k \lambda^k=+\infty,\:\:
\sum_{k=1}^{\infty}(\epsilon^k)^2<\infty,\:\: \sum_{k=1}^{\infty}
\epsilon^k(\lambda^k)^2<\infty
\end{equation}
 then
    the proposed  algorithm will guarantee the following results:
\begin{equation}\label{eq:gradient_summable1_convex}
\begin{aligned}
&\lim_{t\rightarrow\infty}\frac{\sum_{k=0}^t
\epsilon^k\lambda^k\mathbb{E}\left[\left\|\nabla
f(\bar{x}^k)\right\|^2\right]}{\sum_{k=0}^t \epsilon^k \lambda^k}=0,\\
&\lim_{t\rightarrow\infty}\frac{\sum_{k=0}^t \epsilon^k\lambda^k
\mathbb{E}\left[\left\|\frac{\sum_{i=1}^m\nabla f_i(x_i^k)
}{m}\right\|^2\right]}{\sum_{k=0}^t \epsilon^k\lambda^k}=0
\end{aligned}
\end{equation}

Moreover, if in addition, $(\epsilon^k)^{\frac{3}{2}}\lambda^k$ is also summable, i.e.,  $\sum_{k=1}^{\infty}(\epsilon^k)^{\frac{3}{2}}\lambda^k<\infty$,
  then the proposed algorithm will guarantee
\begin{equation}\label{eq:convex_convergence}
\lim_{t\rightarrow\infty}\mathbb{E}\left[
f\left(\frac{\sum_{k=0}^{t}\epsilon^k\lambda^k
x_p^k}{\sum_{k=0}^{t}\epsilon^k\lambda^k}\right)\right]=f(x^{\ast})
\end{equation}
 for any $1\leq p \leq m$.  Note that all expectations are taken with respect to the randomness in stochastic gradients  and quantization up until iteration $k-1$.
\end{Theorem 1}}
\begin{proof}
The derivation of the result in (\ref{eq:gradient_summable1_convex}) is the same as Theorem \ref{theorem_convergence_nonconvex}, so we only consider the derivation of the result in (\ref{eq:convex_convergence}). According to (\ref{eq:bar_x_evolution}), we have the distance
between $\bar{x}^k$ and the optimal solution $x^{\ast}$ evolving as
follows
\begin{equation}\label{eq:bar_x_k+1}
\begin{aligned}
&\mathbb{E}\left[\|\bar{x}^{k+1}-x^{\ast}\|^2\right]\\
&=\mathbb{E}\left[\left\|\bar{x}^k-\epsilon^k\lambda^k\frac{\sum_{i=1}^m g_i^k}{m}-x^{\ast}\right\|^2\right]\\
&= \mathbb{E}\left[\|\bar{x}^k-x^{\ast}\|^2\right]+\mathbb{E}\left[\left\|\epsilon^k\lambda^k \frac{\sum_{i=1}^m g_i^k}{m}\right\|^2\right] \\
&\qquad
-2\mathbb{E}\left[\left\langle\bar{x}^k-x^{\ast},\epsilon^k\lambda^k\frac{\sum_{i=1}^mg_i^k}{m}\right\rangle\right]\\
&= \mathbb{E}\left[\|\bar{x}^k-x^{\ast}\|^2\right]+\mathbb{E}\left[\left\|\epsilon^k\lambda^k \frac{\sum_{i=1}^m g_i^k}{m}\right\|^2\right] \\
&\qquad
-2\mathbb{E}\left[\epsilon^k\lambda^k\frac{\sum_{i=1}^m(g_i^k)^T(\bar{x}^k-x^{\ast})}{m} \right]\\
&\leq \mathbb{E}\left[\|\bar{x}^k-x^{\ast}\|^2\right]+(\epsilon^k\lambda^k)^2G^2\\
&\qquad
-2\mathbb{E}\left[\epsilon^k\lambda^k\frac{\sum_{i=1}^m(g_i^k)^T(\bar{x}^k-x^{\ast})}{m}
\right]
\end{aligned}
\end{equation}
where $\langle\cdot\rangle$ denotes inner product. Note that $G$ is
the upper bound of gradients obtained in Lemma
\ref{Lemma:bounded_gradient}.

Using the convexity of $f_i(\cdot)$, we have the following
relationship for each summand of the last term on the right hand
side of (\ref{eq:bar_x_k+1}):
\begin{equation}\label{eq:g_i^k}
\begin{aligned}
& \mathbb{E}\left[ (g_i^k)^T(\bar{x}^k-x^{\ast})\right]\\
& =\mathbb{E}\left[ (g_i^k)^T(x_i^k-x^{\ast}+\bar{x}^k-x_i^k)\right]\\
 &= \mathbb{E}\left[
 (\nabla f_i(x_i^k))^T(x_i^k-x^{\ast}+\bar{x}^k-x_i^k)\right]\\
&= \mathbb{E}\left[ (\nabla f_i(x_i^k))^T(x_i^k-x^{\ast})\right]+  \mathbb{E}\left[(\nabla f_i(x_i^k))^T(\bar{x}^k-x_i^k)\right]\\
&\geq \mathbb{E}\left[f_i(x_i^k)-f_i(x^{\ast})\right]- G\mathbb{E}\left[\|\bar{x}^k-x_i^k\|\right]\\
&
= \mathbb{E}\left[f_i(x_i^k)-f_i(\bar{x}^k)+f_i(\bar{x}^k)-f_i(x^{\ast})\right]\\
&\qquad\qquad- G\mathbb{E}\left[\|\bar{x}^k-x_i^k\|\right]\\
 &
\geq \mathbb{E}\left[f_i(\bar{x}^k)-f_i(x^{\ast})\right]-2
G\mathbb{E}\left[\|\bar{x}^k-x_i^k\|\right]
\end{aligned}
\end{equation}
where the first inequality used the convexity of $f_i$ and the last
inequality used the relationship $f_i(x_i^k)-f_i(\bar{x}^k)\geq
-G\|\bar{x}^k-x_i^k\|$ from Lemma \ref{lemma:function_gradient} in
the Appendix.

Plugging (\ref{eq:g_i^k}) into (\ref{eq:bar_x_k+1}) yields
\[
\begin{aligned}
\mathbb{E}\left[\|\bar{x}^{k+1}-x^{\ast}\|^2\right]&\leq\mathbb{E}\left[\|\bar{x}^k-x^{\ast}\|^2\right]+ (\epsilon^k\lambda^k)^2G^2 \\
&\qquad
-2\epsilon^k\lambda^k\frac{\sum_{i=1}^m \mathbb{E}\left[(f_i(\bar{x}^k)-f_i(x^{\ast}))\right]}{m}\\
&\qquad +4\epsilon^k\lambda^kG\frac{\sum_{i=1}^m
\mathbb{E}\left[\|\bar{x}^k-x_i^k\|\right]}{m}
\end{aligned}
\]
or
\begin{equation}\label{eq:alpha^k}
\begin{aligned}
 &2\epsilon^k\lambda^k\frac{\sum_{i=1}^m \mathbb{E}\left[(f_i(\bar{x}^k)-f_i(x^{\ast}))\right]}{m}\\
 &\qquad\leq\mathbb{E}\left[\|\bar{x}^k-x^{\ast}\|^2\right]-\mathbb{E}\left[\|\bar{x}^{k+1}-x^{\ast}\|^2\right]+\\
&\quad\qquad
(\epsilon^k\lambda^k)^2G^2+4\epsilon^k\lambda^kG\frac{\sum_{i=1}^m
\mathbb{E}\left[\|\bar{x}^k-x_i^k\|\right]}{m}
\end{aligned}
\end{equation}
Using the fact
\[
\frac{\sum_{i=1}^m
\mathbb{E}\left[(f_i(\bar{x}^k)-f_i(x^{\ast}))\right]}{m}=\mathbb{E}\left[
f (\bar{x}^k)-f (x^{\ast}) \right]
\]
we can rewrite (\ref{eq:alpha^k}) as
\begin{equation}\label{eq:final_convex}
\begin{aligned}
 &2\epsilon^k\lambda^k\mathbb{E}\left[
f (\bar{x}^k)-f (x^{\ast}) \right]\\
 &\qquad\leq\mathbb{E}\left[\|\bar{x}^k-x^{\ast}\|^2\right]-\mathbb{E}\left[\|\bar{x}^{k+1}-x^{\ast}\|^2\right]+\\
&\quad\qquad
(\epsilon^k\lambda^k)^2G^2+4\epsilon^k\lambda^kG\frac{\sum_{i=1}^m
\mathbb{E}\left[\|\bar{x}^k-x_i^k\|\right]}{m}
\end{aligned}
\end{equation}

Summing (\ref{eq:final_convex}) from $k=0$ to $k=t$ yields
\begin{equation}\label{eq:sum}
\begin{aligned}
&2\sum_{k=0}^t\epsilon^k\lambda^k\mathbb{E}\left[(f(\bar{x}^k)-f(x^{\ast}))\right]\\
&\leq
 \mathbb{E}\left[(\bar{x}^0-x^{\ast})^2\right]-\mathbb{E}\left[(\bar{x}^{t+1}-x^{\ast})^2\right]\\
&\qquad + G^2\sum_{k=0}^t (\epsilon^k\lambda^k)^2 \\
&\qquad
+4G\frac{\sum_{k=0}^t\sum_{i=1}^m\epsilon^k\lambda^k\mathbb{E}\left[\|\bar{x}^k-x_i^k\|\right]}{m}
\end{aligned}
\end{equation}
Given that $f(\cdot)$ is a convex function, we always have
\[
f\left(\frac{\sum_{k=0}^t\epsilon^k\lambda^k\bar{x}^k
}{\sum_{k=0}^t\epsilon^k\lambda^k}\right)\leq
\sum_{k=0}^t\frac{\epsilon^k\lambda^kf(\bar{x}^k)}{\sum_{k=0}^t\epsilon^k\lambda^k}
\]
which, in combination with (\ref{eq:sum}), implies
\begin{equation}\label{eq:sum_final}
\begin{aligned}
&
\mathbb{E}\left[f\left(\frac{\sum_{k=0}^t\epsilon^k\lambda^k\bar{x}^k
}{\sum_{k=0}^t\epsilon^k\lambda^k}\right)-f(x^{\ast})\right]\\
&\leq
  \frac{\mathbb{E}\left[\|\bar{x}^0-x^{\ast}\|^2\right]}{2m\sum_{k=0}^t\epsilon^k\lambda^k}-\frac{\mathbb{E}\left[\|\bar{x}^{t+1}-x^{\ast}\|^2\right]}{2m\sum_{k=0}^t\epsilon^k\lambda^k}\\
&\qquad +\frac{G^2\sum_{k=0}^t
(\epsilon^k\lambda^k)^2}{2m\sum_{k=0}^t\epsilon^k\lambda^k}\\
&\qquad
+2G\frac{\sum_{k=0}^t\sum_{i=1}^m\epsilon^k\lambda^k\mathbb{E}\left[\|\bar{x}^k-x_i^k\|\right]}{m\sum_{k=0}^t\epsilon^k\lambda^k}
\end{aligned}
\end{equation}

 Next, we proceed to show that the right hand side of
(\ref{eq:sum_final}) will converge to zero. Based on Lemma
\ref{lemma:lemma_boundedx^2}, we know that
$\mathbb{E}\left[\|\bar{x}^k-x^{\ast}\|\right]$ and
$\mathbb{E}\left[\|\bar{x}^{t+1}-x^{\ast}\|^2\right]$ are bounded,
so the first two terms on the right hand side of
(\ref{eq:sum_final}) will converge to zero under the assumption that
$\epsilon^k\lambda^k$ is  not summable.  The assumption on summable
$(\epsilon^k\lambda^k)^2$    guarantees that the third term on the
right hand side of (\ref{eq:sum_final}) will converge to zero.
Finally, according to Lemma \ref{Lemma:convergence_of_x_i_to_bar_x},
$\mathbb{E}\left[\|\bar{x}^k-x_i^k\|\right]$ is of the order of
$\lambda^k$ or $(\epsilon^k)^{\frac{1}{2}}$, so the last term on the
right hand side of (\ref{eq:sum_final}) will also converge to zero
when the sequences $(\epsilon^k)^2$,
$(\epsilon^k)^{\frac{3}{2}}\lambda^k$, and $\epsilon^k(\lambda^k)^2$
are summable.

Further noting that all $x_p^k$ will converge to each other and
hence to $\bar{x}^k$   according to Lemma
\ref{Lemma:convergence_of_x_i_to_bar_x}, we obtain the statement of
Theorem \ref{Theorem:convergence_convex}.
\end{proof}

In fact, if we can specify the convergence rate of $\epsilon^k$ and
$\lambda^k$, we can further obtain the convergence rate of all
agents to the optimal solution:

\begin{Corollary 2}\label{Theorem:convergence_rate_convex}
If the sequences $\epsilon^k$ and $\lambda^k$ are selected in the form of $\lambda^k
=\frac{a_{1}}{(a_3k+1)^{\delta_1}}$ and $\epsilon^k =
\frac{a_2}{(a_3k+1)^{\delta_2}}$ with $a_1$, $a_2$,  and $a_3$ denoting some
positive constants and  positive exponents $\delta_1$ and $\delta_2$
satisfying $\delta_1+\delta_2\leq1$, $\delta_2>0.5$, and
$2\delta_1+\delta_2>1$,  then all conditions in
(\ref{eq:conditions_convex}) are satisfied and the proposed algorithm will
guarantee (\ref{eq:gradient_summable1_convex}) under Assumptions
\ref{Assumption:gradient}, \ref{assumption:coupling}, and
\ref{Assmption:quantization}. More specifically, the
convergence rate of gradients satisfies
\begin{equation}\label{eq:rate_in_nonconvex}
\begin{aligned}
&\frac{\sum_{k=0}^{t}\epsilon^k \lambda^k\mathbb{E}\left[
\left\|\nabla
f(\bar{x}^k)\right\|^2\right]}{\sum_{k=0}^{t}\epsilon^k \lambda^k}
\\
&\qquad+\frac{\sum_{k=0}^{t}\epsilon^k  \lambda^k\mathbb{E}\left[
\left\|\frac{\sum_{i=1}^m\nabla f_i(x_i^k)
}{m}\right\|^2\right]}{\sum_{k=0}^{t}\epsilon^k\lambda^k}\\
&\qquad+
\frac{2\left(\mathbb{E}\left[f(\bar{x}^{t+1})\right]-\mathbb{E}\left[f(\bar{x}^{0})\right]\right)}{\sum_{k=0}^{t}\epsilon^k
\lambda^k}\\
& =\mathcal{O}\left(\frac{1}{(t+1)^{\delta}}\right)
\end{aligned}
\end{equation}
where $\delta=\min\{2\delta_1,\delta_2\}$.

If in addition,  $\delta_1$ and $\delta_2$
satisfy
$\delta_1+\frac{3}{2}\delta_2\geq 1$,
then  the convergence rate of function values satisfies
\begin{equation}\label{eq:convergence_rate_convex}
\begin{aligned}
&\mathbb{E}\left[ f\left(\frac{\sum_{k=0}^{k}\epsilon^k\lambda^k
x_p^k}{\sum_{k=0}^{t}\epsilon^k\lambda^k}\right)-f(x^{\ast})\right]\\
&+\frac{\mathbb{E}\left[\|\bar{x}^{t+1}-x^{\ast}\|^2\right]-\mathbb{E}\left[\|\bar{x}^0-x^{\ast}\|^2\right]}{2m\sum_{k=0}^t\epsilon^k\lambda^k}
=\mathcal{O}\left(\frac{1}{(t+1)^{\delta}}\right)
\end{aligned}
\end{equation}
where $\delta=\min\{ \delta_1, \frac{1}{2}\delta_2\}$ for any $1\leq
p\leq m$. { Note that all expectations are taken with respect to the randomness in stochastic gradients  and quantization up until iteration $k-1$.}
\end{Corollary 2}
\begin{proof}
The statement for the convergence rate of  gradients follows Corollary \ref{Theorem:convergence_rate}.
To arrive at the statement on the convergence rate of the function value, one can follow  the line of derivation in the proof of Theorem
\ref{Theorem:convergence_convex}. More specifically, under the
conditions of Theorem \ref{Theorem:convergence_rate_convex}, we can
obtain that in (\ref{eq:sum_final}), the numerators of the second
and third terms on the right hand side will decay with a rate of no
less than $\mathcal{O}\left(\frac{1}{(t+1)^{\delta}}\right)$ with
$\delta=\min\{2(\delta_1+\delta_2),2\delta_1+\delta_2,\delta_1+\frac{3}{2}\delta_2\}$.
We further note that the denominator $\epsilon^k\lambda^k$   decays
with the rate of $\delta_1+\delta_2$, and hence that
 the left hand side of (\ref{eq:convergence_rate_convex})  decays with a rate of $\delta=\min\{
\delta_1,\frac{1}{2}\delta_2\}$ as  in the statement of  the
theorem.
\end{proof}

\section{Privacy Analysis}

In this section, we show that  our algorithm's  robustness to
aggressive quantization effects can   be leveraged to enable
rigorous differential privacy. More specifically, under a ternary
quantization scheme   which quantizes any value to three numerical
levels, we will prove that  our decentralized optimization algorithm
can enable rigorous differential privacy without losing provable
convergence accuracy. To the best of our knowledge, this is the
first time both strict differential privacy and provable convergence
accuracy are achieved in decentralized stochastic optimization.

The ternary quantization scheme is defined as follows:
\begin{Definition 1}
The ternary quantization scheme quantizes a vector
$x=[x_1,x_2,\cdots,x_d]^T\in \mathbb{R}^d$ as follows
\[
\mathcal{Q}(x)=\left[q_1,\,q_2,\,\cdots,q_d\right],\quad q_i= r{\rm sign}(x_i)b_i,\quad \forall 1\leq i\leq d
\]
where $r$ is a design parameter no less than the $\ell_\infty$ norm
$\|x\|_{\infty}$ of $x$, ${\rm sign}$ represents the sign of a
value, and $b_i$ ($1\leq i\leq d$) are independent binary variables
following the Bernoulli distribution
\[
\left\{\begin{array}{ccc}P(b_i=1|x)&=&|x_i|/r\\P(b_i=0|x)&=&1-|x_i|/r\end{array}\right.
\]
with $P(\cdot)$ denoting the probability distribution.
\end{Definition 1}

Such ternary quantization has been applied in distributed stochastic
optimization, in, e.g.,
\cite{reisizadeh2019exact,alistarh2017qsgd,wen2017terngrad}.
However,   none of these results use quantization effects to achieve
strict differential privacy. Now we show that using the ternary
quantization, our decentralized stochastic optimization algorithm
can achieve $(0,\delta)$-differential privacy while maintaining
provable convergence accuracy:
\begin{Theorem 1}\label{th:privacy}
Under  Assumptions  \ref{Assumption:gradient},\ref{assumption:coupling} in the non-convex case, or Assumptions \ref{Assumption:gradient},\ref{assumption:coupling},\ref{Assumtion:convex} in the convex case,  the  ternary quantization scheme defined
in Definition 2  { achieves $(0,\frac{1}{r})$-differential privacy for
individual agents' gradients in every iteration while ensuring convergence.}
\end{Theorem 1}
\begin{proof}
It can be easily verified that the ternary quantization scheme
satisfies the conditions in Assumption \ref{Assmption:quantization}.
So the decentralized optimization algorithm will have provable
convergence accuracy according to Theorem
\ref{theorem_convergence_nonconvex} and Theorem
\ref{Theorem:convergence_convex},
 and we only need to prove that $(0,\frac{1}{r})$-differential privacy can be obtained for individual agents' gradients under such a
quantization scheme.

From the proposed algorithm in (\ref{eq:proposed_algorithm}), it can
be seen that for an individual agent $i$, its gradient $g_i^k$ can
be viewed as a function of  all variables $x_i^k$ ($1\leq i\leq m$).
Therefore, using differential privacy's robustness to
post-processing operations \cite{dwork2014algorithmic}, if we can
prove that the ternary quantization scheme can enable
$(0,\frac{1}{r})$-differential privacy for $x_i^k$, then we have
that the ternary quantization scheme can enable
$(0,\frac{1}{r})$-differential privacy for individual agents'
gradients.

{According to the mechanism of ternary quantization, it can be obtained that depending on the sign of $x_i^k$, the quantized value can have different distributions:
\[
\begin{aligned}
\left\{\begin{array}{lcc}P(q_i=r|x)&=&|x_i|/r\\P(q_i=0|x)&=&1-|x_i|/r\\P(q_i=-r|x)&=&0\end{array}\right.\quad \textrm{when} \:\:x_i\geq 0
\end{aligned}
\]
and
\[
\begin{aligned}
\left\{\begin{array}{lcc}P(q_i=r|x)&=&0\\P(q_i=0|x)&=&1-|x_i|/r\\P(q_i=-r|x)&=&|x_i|/r\end{array}\right.\quad \textrm{when} \:\:x_i< 0
\end{aligned}
\]
{Furthermore,  given that the quantization of one element is independent of that of  other elements, i.e., the quantization errors for different elements are independent of each other, we can consider the per-step privacy of different elements of  $x$  separately.}
Therefore, according to Definition 1, to prove that $(0,\frac{1}{r})$-differential privacy is achieved, i.e., $\left|P(q_i\in S|y_i)-P(q_i \in S|x_i)\right|\leq \frac{1}{r}$ for all $S\in \{r,\,0,-r\}$ and all $x$, $y$ with $\|x-y\|_1\leq 1$, we divide the derivation into two cases: 1) $x_i$ and $y_i$ are of the same sign, i.e., both $x_i$ and $y_i$ are nonnegative or both $x_i$ and $y_i$ are negative; 2) $x_i$ and $y_i$ are of different signs, i.e., either $x_i\geq 0$, $y_i<0$ is true or $x_i< 0$, $y_i\geq 0$  is true.

{\bf Case 1}: $x_i$ and $y_i$ are of the same sign, i.e.,  both $x_i$ and $y_i$ are nonnegative or both $x_i$ and $y_i$ are negative. Without loss of generality, we assume that both $x_i$ and $y_i$ are nonnegative. It can be easily verified that the same result can be obtained if both $x_i$ and $y_i$ are negative.

 Based
on the mechanism of ternary quantization, it can be obtained that
\[
\begin{aligned}
&\sup_{\|x-y\|_1\leq 1}\big|
P(q_i=r|x)-P(q_i=r|y)\big|\\
&\qquad =\sup_{\|x-y\|_1\leq 1}\left|\frac{|x_i|-|y_i|}{r}\right|\leq \frac{\|x-y\|_1}{r}\leq \frac{1}{r},\\
&\sup_{\|x-y\|_1\leq 1} \left|P(q_i=0|x)-P(q_i=0|y)\right|\\
&\qquad =\sup_{\|x-y\|_1\leq 1}\left|\frac{(r-|x_i|)-(r-|y_i|)}{r}\right|\leq
\frac{\|x-y\|_1}{r}\leq \frac{1}{r},\\
&\sup_{\|x-y\|_1\leq 1}\big|
P(q_i=-r|x)-P(q_i=-r|y)\big|\\
&\qquad =\sup_{\|x-y\|_1\leq 1}|0-0|\leq \frac{1}{r}
\end{aligned}
\]

In a similar way, one can obtain the same relationship when both $x$ and $y$ are negative.

{\bf Case 2}: $x_i$ and $y_i$ are of different signs, i.e., either $x_i\geq 0$, $y_i<0$ is true or $x_i< 0$, $y_i\geq 0$  is true. Without loss of generality, we assume that $x_i\geq 0$, $y_i<0$ is true.  It can be easily verified that the same result can be obtained if $x_i< 0$, $y_i\geq 0$ is true.

Under the constraint $x_i\geq 0$ and $y_i<0$, it can be obtained that  $|x_i|\leq 1$ and $|y_i|\leq 1$ must hold for all $x$ and $y$ satisfying $\|x-y\|_1\leq 1$. Therefore, based
on the mechanism of ternary quantization, it can be obtained that
\[
\begin{aligned}
&\sup_{\|x-y\|_1\leq 1}\big|
P(q_i=r|x)-P(q_i=r|y)\big|\\
&\qquad =\sup_{\|x-y\|_1\leq 1}\left|\frac{|x_i|}{r}-0\right|\leq \frac{|x_i|}{r}\leq \frac{1}{r},\\
&\sup_{\|x-y\|_1\leq 1} \big|P(q_i=0|x)-P(q_i=0|y)\big|\\
&\qquad =\sup_{\|x-y\|_1\leq 1}\left|\frac{(r-|x_i|)-(r-|y_i|)}{r}\right|\leq
\frac{\|x-y\|_1}{r}\leq \frac{1}{r},\\
&\sup_{\|x-y\|_1\leq 1}\big|
P(q_i=-r|x)-P(q_i=-r|y)\big|\\
&\qquad =\sup_{\|x-y\|_1\leq 1}\left|0-\frac{|y_i|}{r}\right|\leq \frac{|y_i|}{r}\leq \frac{1}{r}
\end{aligned}
\]
Summarizing the results in Case 1 and Case 2,   we always have $(0, \frac{1}{r})$-differential privacy for the quantization
 input $x_i^k$ for every individual agent $i$. Further using the robustness of differential privacy  to
post-processing operations \cite{dwork2014algorithmic} yields that
we have $(0, \frac{1}{r})$-differential privacy for all agents'
gradients.}
\end{proof}

 Since in
$(\epsilon, \delta)$-differential privacy, the strength of privacy
protection increases with a decrease in $\epsilon$ and $\delta$, the
achieved $(0, \delta)$-differential privacy is stronger than
commonly used $(\epsilon, \delta)$-differential privacy.
Furthermore,  we can see that a larger threshold value $r$ reduces
$\frac{1}{r}$, and hence will lead to a stronger privacy protection.
This is intuitive as a larger $r$ will mean a higher  probability of
no transmission under a given input (since the value to be
transmitted
 is 0). {However, a larger $r$ will also slow  down convergence, as
illustrated in Fig. \ref{fig:comp_convex} in the numerical
simulation section.}

{
\begin{Remark 1}
From the derivation, it can be verified that the same $(0,\frac{1}{r})$-differential privacy can still be obtained when the $\ell_1$ norm in Definition 1 is replaced with any $\ell_p$ norm defined by $\|x\|_p=\left(|x_1|^p+|x_2|^p+\cdots+|x_m|^p\right)^{1/p}$ with $p\geq 1$.
\end{Remark 1}

\begin{Remark 1}
From the derivation, it can also be seen that the stochastic nature of the quantizer is crucial for enabling differential privacy on shared messages.
\end{Remark 1}
}

\begin{Remark 1}
Note that the proposed algorithm can guarantee the privacy of all
participating agents even when an adversary has access to all shared
messages in the network. This is in distinct difference from
existing accuracy-friendly privacy  solutions  (in, e.g.,
\cite{li2020privacy,yan2012distributed,lou2017privacy,gade2018privacy}
for decentralized deterministic convex optimization) that will fail
to protect privacy   when  an adversary has  access to all shared
messages in the network.
\end{Remark 1}

{
\begin{Remark 1}
   Note that an adversary can obtain the information that the quantizer input is no larger than $r$.
\end{Remark 1}
}

{
\begin{Remark 1}
   Theorem \ref{th:privacy} provides privacy guarantee for one quantization operation, i.e., one iteration. The cumulative privacy loss (budget) increases roughly at a rate of $\sqrt{T}$ for $T$ iterations, according to the composition theorem for differential privacy \cite{kairouz2015composition}.
\end{Remark 1}
}

{
\begin{Remark 1}
 The proposed results are significantly different from \cite{agarwal2018cpsgd}. First, we consider the fully decentralized scenario with no servers, whereas \cite{agarwal2018cpsgd} addresses the scenario with a server-client architecture, whose convergence analysis is fundamentally  different from the server-free decentralized case. Moreover, the privacy mechanism in \cite{agarwal2018cpsgd} still falls within the conventional noise-injecting framework for differential privacy since it considers quantization and privacy separately (\cite{agarwal2018cpsgd} uses a dedicated noise mechanism to generate noise and then injects the noise on the quantization output, although  binomial noise is used instead of commonly used Gaussian noise), whereas the  approach in this paper  exploits the quantization error directly to achieve privacy   and hence avoids any dedicated noise-injection mechanism.
\end{Remark 1}
}

Under the ternary quantization scheme,  any transmitted value is
represented as a ternary vector with three possible values
$\{-r,0,r\}$. So to transmit a value,  instead of transmitting
32-bits, which is the typical number of bits to represent a value in
modern computing devices, we could instead only transmit much fewer
bits in addition to the threshold value. So theoretically ternary
quantization can reduce the traffic by a factor of
$\frac{32}{\log_2(3)}=20.18\times$. Therefore, our decentralized
optimization algorithm with ternary quantization  can have
communication efficiency, strict $(0, \delta)$-differential privacy,
as well as provable convergence accuracy simultaneously. To the best
of our knowledge, this is the first decentralized optimization
algorithm able to achieve these three goals simultaneously.



\section{Numerical Experiments}

In this section, we  evaluate the performance of our algorithm using
numerical experiments. We will consider both the convex
objective-function case and the non-convex objective-function case.
\subsection{Convex case}
For the case of convex objective functions, we consider a canonical
decentralized estimation problem where a sensor network of $m$
sensors collectively estimate an unknown parameter
$\theta\in\mathbb{R}^d$, which can be formulated as an empirical
risk minimization problem. More specifically, we assume that each
sensor $i$ has $n_i$ noisy measurements of the parameter
$z_{ij}=M_i\theta+w_{ij}$ for $j=\{1,2,\cdots,n_i\}$ where
$M_i\in\mathbb{R}^{s\times d}$ is the measurement matrix of agent
$i$ and $w_{ij}$ is measurement noise associated with measurement
$z_{ij}$. Then the estimation of the parameter $\theta$ can be
solved using the decentralized optimization problem formulated in
(\ref{eq:optimization_formulation1}), with each $f_i(\theta)$ given
by
\[
f_i(\theta)=\frac{1}{n_i}\sum_{j=1}^{n_i}\|z_{ij}-M_i\theta\|^2+r_i\|\theta\|^2
\]
where $r_i$ is a non-negative regularization parameter.

\begin{figure}
    \begin{center}
        \includegraphics[width=0.3\textwidth]{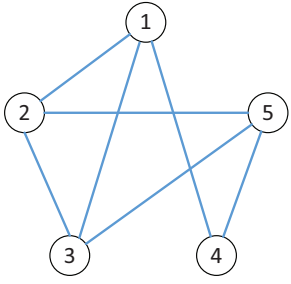}
    \end{center}
    \caption{The interaction topology of the network.}
    \label{fig:topology}
\end{figure}

\begin{figure}
    \begin{center}
        \includegraphics[width=0.45\textwidth]{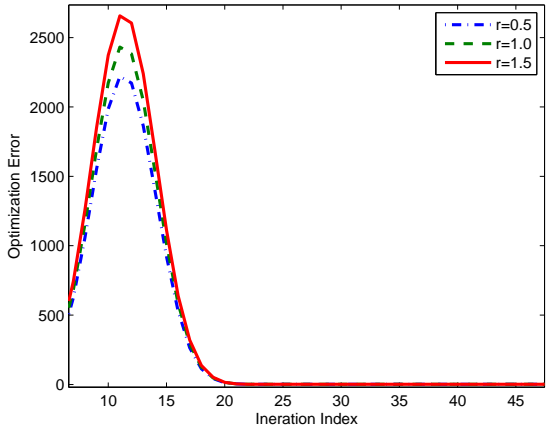}
    \end{center}
    \caption{Comparison of convergence performance under different thresholds of the quantization scheme. Here the optimization error is defined as $\|x^{\ast}-x^k\|$.}
    \label{fig:comp_convex}
\end{figure}

We assume that the network  consists  of five agents interacting on
a graph depicted in Fig. \ref{fig:topology}. The dimension $s$ was
set to 3 and the dimension $d$ was set to 2. $n_i$ was set to 100
for all $i$. $w_{ij}$ were assumed to be uniformly distributed in
$[0,\,1]$. To evaluate the performance of our proposed decentralized
stochastic optimization algorithm, we set
$\lambda^k=\frac{1}{(0.3k+1)^{0.3}}$ and
$\epsilon^k=\frac{1}{(0.3k+1)^{0.6}}$. It can be verified that the
parameters  satisfy the  conditions required in Theorem
\ref{Theorem:convergence_convex} and Corollary
\ref{Theorem:convergence_rate_convex}.  The evolution of the
estimation error averaged over 100 runs is illustrated   in Fig.
\ref{fig:comp_convex}, where we show the results under three
different threshold values of the quantization scheme. It can be
seen that a larger threshold tends to bring a larger overshoot in
the optimization process.

\subsection{Non-convex case}
\begin{figure}
    \begin{center}
        \includegraphics[width=0.5\textwidth]{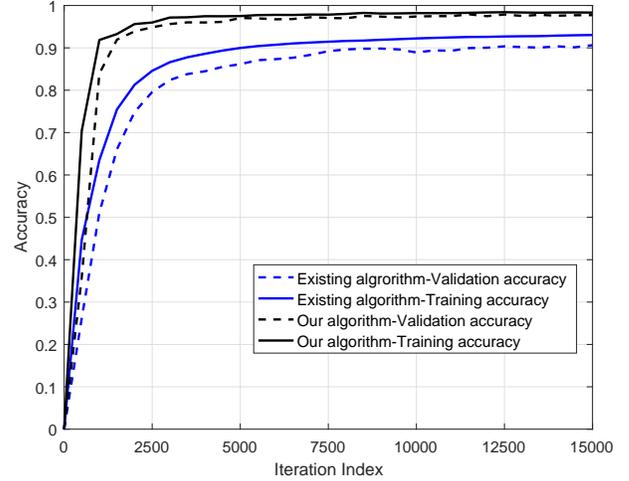}
    \end{center}
    \caption{Comparison of CNN training/validation performance between our algorithm and the conventional decentralized stochastic optimization algorithm in \cite{lian2017can}. }
    \label{fig:nonconvex}
\end{figure}

We use the decentralized training of a convolutional neural network
(CNN)  to evaluate the performance of our proposed decentralized
stochastic optimization algorithm in non-convex optimization. More
specially, we consider five agents interacting on a topology
depicted in Fig. \ref{fig:topology}. The agents collaboratively
train a CNN using the MNIST data set \cite{MNIST}, which is a large
benchmark database of handwritten digits widely used for training
and testing in the field of machine learning \cite{deng2012mnist}.
Each agent has a local copy of the CNN. The CNN has 2 convolutional
layers with 32 filters, and then two more convolutional layers with
64 filters each followed by a dense layer with 512 units. Each agent
has access to a portion of the MNIST data set, which was further
divided into two subsets for training and validation, respectively.
We set the optimization parameters   as $\lambda^k=
\frac{1}{(0.001k+1)^{0.3}}$ and
$\epsilon^k=\frac{1}{(0.001k+1)^{0.7}}$.  For the adopted CNN model,
the dimension of gradient, $d$, is equal to $1,676,266$.  It can be
verified that the parameters satisfy the  conditions required in
Theorem \ref{theorem_convergence_nonconvex} and Corollary
\ref{Theorem:convergence_rate}.  The evolution of the training and
validation accuracies averaged over 100 runs are illustrated by the
solid and dashed black lines in Fig. \ref{fig:nonconvex}. To compare
the convergence performance of our algorithm with the conventional
decentralized stochastic optimization algorithm, we also implemented
the decentralized stochastic optimization algorithm in
\cite{lian2017can} to train the same CNN under the same quantization
scheme, whose average training and validation accuracies over 100
runs are represented by the solid and dashed blue lines in Fig.
\ref{fig:nonconvex}. It can be seen that the proposed algorithm has
a faster  converging rate as well as  better training/validation
accuracy in the presence of quantization effects.

\begin{figure}
    \begin{center}
        \includegraphics[width=0.5\textwidth]{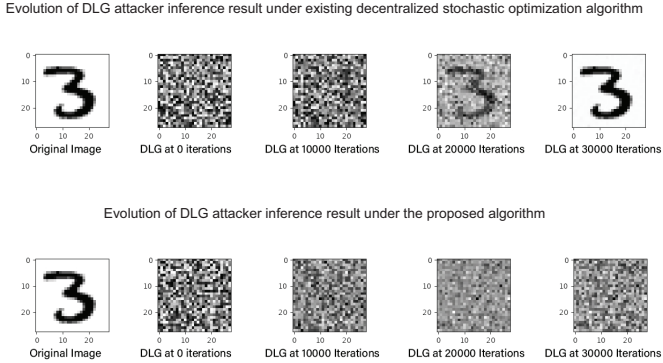}
    \end{center}
    \caption{Comparison of DLG attacher's inference results under existing decentralized stochastic optimization algorithm in \cite{lian2017can} and our algorithm.}
    \label{fig:DLG_image}
\end{figure}

\begin{figure}
    \begin{center}
        \includegraphics[width=0.45\textwidth]{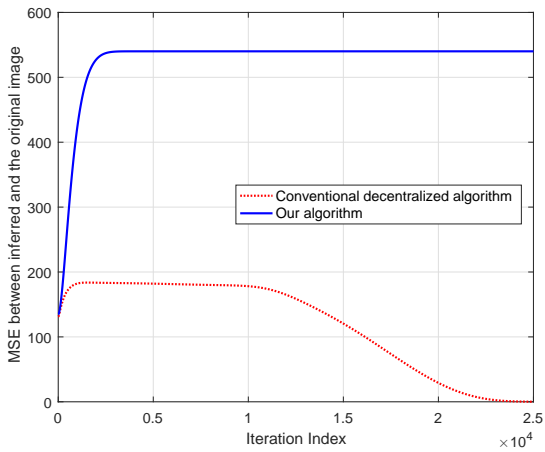}
    \end{center}
    \caption{Comparison of DLG attacher's inference errors under existing decentralized stochastic optimization algorithm in \cite{lian2017can} and our algorithm.}
    \label{fig:DLG_MSE}
\end{figure}

To show that the proposed algorithm can indeed protect the privacy
of participating agents, we also implemented a privacy attacker
which tries to infer the raw image of participating agents using
received information. The attacker implements the DLG attack model
 proposed in \cite{zhu2019deep}, which is the most powerful
inference algorithm reported to date in terms of reconstructing
exact raw data from shared gradients/model updates. The attacker was
assumed to be able to eavesdrop all messages shared among the
agents. Fig. \ref{fig:DLG_image} shows that the attacker could
effectively recover the original training image from shared model
updates in the conventional stochastic optimization algorithm in
\cite{lian2017can} {that does not take privacy protection into consideration}. However, under the proposed algorithm and
quantization effects, the attacher failed to infer the original
training image through information shared in the network. This is
also corroborated by the attacker's inference performance measured
by the mean-square error (MSE) between the inference result  and the
original image. More specifically, as illustrated in Fig.
\ref{fig:DLG_MSE}, the attacker eventually inferred the raw image
accurately as its estimation error converged to zero. However, the
proposed approach successfully thwarted the attacker as attacker's
estimation error was always large.

\section{Conclusions}
The paper has presented a decentralized stochastic optimization
algorithm that is robust  to aggressive quantization effects, which
 enables the exploitation of aggressive quantization effects to
obfuscate shared information and hence enables privacy protection in
decentralized stochastic optimization without losing provable
convergence accuracy. Based on this result, this paper, for the
first time, proposes and achieves ternary-quantization based
rigorous ($0,\delta$)-differential privacy without losing provable
convergence accuracy in decentralized stochastic optimization.  The
results are applicable in both the convex optimization case and the
non-convex optimization case. The ternary quantization scheme also
leads to significant reduction in communication overhead. Our
approach appears to be the first to achieve rigorous differential
privacy, communication efficiency, and provable convergence accuracy
simultaneously in decentralized stochastic optimization. Both
simulation results for a convex decentralized optimization problem
and numerical experimental results for machine learning on a
benchmark image dataset confirm the effectiveness of the proposed
approach.

The paper assumes smooth gradients and does not consider potential
constraints between optimization variables, as, for example, in
\cite{cao2020decentralized}. In the future, we plan to extend the
results to more general non-smooth and constrained decentralized
optimization problems.

\section*{Acknowledgement}
The authors would like to thanks Ben Liggett for the help in numerical experiments. They would also like to thank the anonymous reviewers, whose comments helped improve  the paper.

\section*{Appendix}

\subsection{Some preliminary results}\label{Appendix:preliminary}
\begin{Lemma
1}\label{Lemma:convergence_lemma}\cite{robbins1971convergence} Let
$\{v^k\}$ be a non-negative sequence  satisfying the following
relationship for all $k\geq 0$:
\begin{equation}\label{eq:sequence_v^k1}
v^{k+1}\leq (1+a^k)v^k+w^k
\end{equation}
where sequences $a^k\geq 0$ and $w^k\geq 0$ satisfy
$\sum_{k=0}^{\infty}a^k<\infty$ and $\sum_{k=0}^{\infty}w^k<\infty$,
respectively. Then the sequence $\{v^k\}$ will converge to a finite
value $v\geq 0$.
\end{Lemma 1}

\begin{Lemma
2}\label{Lemma:decay_rates}\cite{kar2013distributed,george2019distributed}
Let $\{v^k\}$ be a non-negative sequence which satisfies the
following relationship for all $k\geq 0$:
\begin{equation}\label{eq:sequence_v^k2}
v^{k+1}\leq (1-r_1^k)v^k+r_2^k
\end{equation}
with sequences $r_1^k\geq 0$ and $r_2^k\geq 0$   satisfying
\[
\frac{C_1}{(C_3k+1)^{ \gamma_1}}\leq r_1^k\leq 1,\quad
\frac{C_2}{(C_3k+1)^{ \gamma_2}}\leq r_2^k\leq 1
\]
for some $C_1>0$, $C_2>0$, $C_3>0$,  $0\leq \gamma_1<1$, and
$\gamma_1<\gamma_2$. Then $\lim_{k\rightarrow
\infty}(k+1)^{\gamma_0}v^k=0$ holds for all $0\leq
\gamma_0<\gamma_2-\gamma_1$.
\end{Lemma 2}

\begin{Lemma 2}\label{lemma:function_gradient}\cite{nedic2018network}
Suppose $h:\mathbb{R}^d\rightarrow\mathbb{R}$ is a convex function
with  gradient bounded by $G$. Then we have
\[
|h(y)-h(x)|\leq G\|y-x\|
\]
for any $x,y\in\mathbb{R}^d$
\end{Lemma 2}

\subsection{\bf Proof of Lemma \ref{lemma:lemma_boundedx^2}}\label{Appendix:Lemma_bounded}
According to Lemma \ref{Lemma:convergence_lemma} in the Appendix, to
prove that $\mathbb{E}\left[\|x^k\|^2\right]$ is bounded, we only
need to prove that under the conditions in Lemma
\ref{lemma:lemma_boundedx^2}, it satisfies the inequality in
(\ref{eq:sequence_v^k1}) in the Appendix.

For the convenience of analysis, we first define the augmented
versions of $x^{\ast}$ and $\bar{x}^k$:
\begin{equation}\label{eq:augmented_barx_Omega}
\hat{x}^{\ast}\triangleq  {\bf 1}_m\otimes x^{\ast},\quad
\hat{\bar{x}}^k\triangleq  {\bf 1}_m\otimes\bar{x}^k
\end{equation}
where ${\bf 1}_m$ denotes an $m$ dimensional column vector with all
entries equal to 1.

  Using the
inequality $(x+y)^2\leq 2x^2+2y^2$, which holds for any
$x,y\in\mathbb{R}$, we can obtain
\begin{equation}\label{eq:x_bar_x_tilde_x}
\begin{aligned}
\|x^k\|^2&=\|\hat{\bar{x}}^k-\hat{x}^{\ast}+x^k-\hat{\bar{x}}^k+\hat{x}^{\ast}\|^2 \\
&\leq (\|\hat{\bar{x}}^k-\hat{x}^{\ast}+x^k-\hat{\bar{x}}^k\|+\|\hat{x}^{\ast}\|)^2\\
&\leq 2\|\hat{\bar{x}}^k-\hat{x}^{\ast}+x^k-\hat{\bar{x}}^k\|^2+2\|\hat{x}^{\ast}\|^2\\
& \leq 2(\|\hat{\bar{x}}^k-\hat{x}^{\ast}\|+\|x^k-\hat{\bar{x}}^k\|)^2+2\|\hat{x}^{\ast}\|^2\\
&\leq 4 \|\hat{\bar{x}}^k-\hat{x}^{\ast}\|^2+4\|x^k-\hat{\bar{x}}^k\|^2+2\|\hat{x}^{\ast}\|^2\\
\end{aligned}
\end{equation}

Because $\hat{x}^{\ast}$ is a constant, we will prove the
boundedness of $\mathbb{E}\left[\|x^k\|^2\right]$ by proving that
$\mathbb{E}\left[ \|\hat{\bar{x}}^k-\hat{x}^{\ast}\|^2+
\|x^k-\hat{\bar{x}}^k\|^2 \right]$ is bounded. Our derivation will
follow three steps: in Step I and Step II, we  study the respective
evolution of  $\mathbb{E}\left[
\|\hat{\bar{x}}^k-\hat{x}^{\ast}\|^2\right]$ and
$\mathbb{E}\left[\|x^k-\hat{\bar{x}}^k\|^2 \right]$ under our
proposed algorithm in (\ref{eq:proposed_algorithm}); in Step III, we
show that $\mathbb{E}\left[ \|\hat{\bar{x}}^k-\hat{x}^{\ast}\|^2+
\|x^k-\hat{\bar{x}}^k\|^2 \right]$ is bounded by combining the
relationship obtained in Step I and Step II.

{\bf Step I}: We first consider $\mathbb{E}\left[
\|\hat{\bar{x}}^k-\hat{x}^{\ast}\|^2\right]$, which is equal to
$m\mathbb{E}\left[ \| {\bar{x}}^k- x^{\ast}\|^2\right]$ according to
the definition in (\ref{eq:augmented_barx_Omega}).

From (\ref{eq:bar_x_evolution}), we have
\begin{equation}\label{eq:bar_x_inequality}
\begin{aligned}
\|\bar{x}^{k+1}-x^{\ast}\|^2=&
\left\|\bar{x}^k-x^{\ast}-\epsilon^k\lambda^k\frac{\sum_{i=1}^m g_i^k}{m}\right\|^2\\
 \leq&
\left(\left\|\bar{x}^k-x^{\ast}\right\|+\left\|\epsilon^k\lambda^k\frac{\sum_{i=1}^m
g_i^k}{m}\right\|\right)^2
\end{aligned}
\end{equation}

Using the inequality $(x+y)^2\leq (1+\nu)x^2+(1+\frac{1}{\nu})y^2$,
which holds for any $x,y\in\mathbb{R}$ and $\nu>0$, we can obtain
the following relationship from (\ref{eq:bar_x_inequality}) by
setting $\nu$ to $(\epsilon^k)^2$
\begin{equation}\label{eq:after_applying_first}
\begin{aligned}
\|\bar{x}^{k+1}-x^{\ast}\|^2&\leq
 \left(1+(\epsilon^k)^2\right)\left\|\bar{x}^k-x^{\ast}\right\|^2\\
 &\quad+\left(1+\frac{1}{(\epsilon^k)^2}\right)\left\|\epsilon^k\lambda^k\frac{\sum_{i=1}^m
g_i^k}{m}\right\|^2\\
&=\left(1+(\epsilon^k)^2\right)\left\|\bar{x}^k-x^{\ast}\right\|^2\\
&\quad+\left((\epsilon^k\lambda^k)^2+ (\lambda^k)^2
\right)\left\|\frac{\sum_{i=1}^m
g_i^k}{m}\right\|^2\\
&\leq\left(1+(\epsilon^k)^2\right)\left\|\bar{x}^k-x^{\ast}\right\|^2\\
&\quad+\left((\epsilon^k\lambda^k)^2+
(\lambda^k)^2 \right)G^2\\
 \end{aligned}
\end{equation}
where we used the result that the gradient is bounded by $G$ from
Lemma \ref{Lemma:bounded_gradient}.

{\bf Step II:} We next consider
$\mathbb{E}\left[\|x^k-\hat{\bar{x}}^k\|^2 \right]$. From
(\ref{eq:proposed_algorithm_vector_form}) and
(\ref{eq:bar_x_evolution}), we can obtain the dynamics of
$x^k-\hat{\bar{x}}^k$ based on the fact $ A^k\bar{x}^k=\bar{x}^k$:
\begin{equation}
\begin{aligned}
x^{k+1}-\hat{\bar{x}}^{k+1}&= ({A}^k\otimes I_d)
(x^k-\hat{\bar{x}}^k)-\epsilon^k\lambda^k  ({M}\otimes I_d)  g^k\\
&\quad + \epsilon^k (L_w\otimes I_d) V^k
\end{aligned}
\end{equation}
where $  M= \left(I-\frac{{\bf 1} {\bf 1}^T}{m}\right) $ and the
other parameters are given in
(\ref{eq:proposed_algorithm_vector_form}). Therefore, we have
\begin{equation}
\begin{aligned}
&\|x^{k+1}-\hat{\bar{x}}^{k+1}\|^2\\
&=\|(A^k\otimes I_d) (x^k-\hat{\bar{x}}^k)-\epsilon^k\lambda^k (M
\otimes I_d)  g^k\|^2
\\
&\quad +\|\epsilon^k (L_w \otimes I_d) V^k\|^2+\\
&\quad 2\left\langle (A^k\otimes I_d)
(x^k-\hat{\bar{x}}^k)-\epsilon^k\lambda^k (M\otimes I_d) g^k ,
\epsilon^k (L_w\otimes I_d) V^k \right\rangle
\end{aligned}
\end{equation}
i.e.,
\begin{equation}\label{eq:consensus_first}
\begin{aligned}
&\mathbb{E}\left[\|x^{k+1}-\hat{\bar{x}}^{k+1}\|^2\right]\\
&=\mathbb{E}\left[\|(A^k\otimes I_d)
(x^k-\hat{\bar{x}}^k)-\epsilon^k\lambda^k (M \otimes
I_d)g^k\|^2\right]\\
&\quad  +\mathbb{E}\left[\|\epsilon^k (L_w\otimes I_d) V^k \|^2\right]
\end{aligned}
\end{equation}
where we used the fact that $V^k$ is uncorrelated noise with
expectation equal to zero.

It can be verified that the following relationship holds
\[
\begin{aligned}
\|(A^k\otimes I_d)& (x^k-\hat{\bar{x}}^k)-\epsilon^k\lambda^k
(M\otimes I_d) g^k\|\\
 &\leq \|(A^k\otimes I_d)
(x^k-\hat{\bar{x}}^k)\|+\|\epsilon^k\lambda^k (M\otimes I_d)
 g^k\|\\
&\leq (1-\epsilon^k\rho)\|
x^k-\hat{\bar{x}}^k\|+\|\epsilon^k\lambda^k (M\otimes I_d)  g^k\|\\
& \leq (1-\epsilon^k\rho)\|
x^k-\hat{\bar{x}}^k\|+\epsilon^k\lambda^k\| g^k\|
\end{aligned}
\]
where the second inequality used  the doubly-stochastic property of
$A^k$ and Lemma 4.4 of \cite{kar2013distributed} with $\rho$ the
second largest eigenvalue of $L_w$, and the third inequality used the
fact $\|M\|=1$. Therefore, we have
\[
\begin{aligned}
\|&(A^k\otimes I_d) (x^k-\hat{\bar{x}}^k)-\epsilon^k\lambda^k
(M\otimes I_d)   g^k\|^2 \\
& \leq (1+\nu)(1-\epsilon^k\rho)^2\|
x^k-\hat{\bar{x}}^k\|^2+(1+\frac{1}{\nu})(\epsilon^k\lambda^k)^2\|
g^k\|^2
\end{aligned}
\]
based on the inequality $(x+y)^2\leq
(1+\nu)x^2+(1+\frac{1}{\nu})y^2$, which holds for any
$x,y\in\mathbb{R}$ and $\nu>0$. Setting $\nu$ as $\epsilon^k\rho$,
we further have
\begin{equation}\label{eq:Consensus_error}
\begin{aligned}
&\|(A^k\otimes I_d) (x^k-\hat{\bar{x}}^k)-\epsilon^k\lambda^k
(M\otimes I_d)  g^k\|^2 \\
& \leq (1+\epsilon^k\rho)(1-\epsilon^k\rho)^2\|
x^k-\hat{\bar{x}}^k\|^2 +
(1+\frac{1}{\epsilon^k\rho})(\epsilon^k\lambda^k)^2\|
g^k\|^2\\
&=(1-(\epsilon^k)^2\rho^2)(1-\epsilon^k\rho)\|
x^k-\hat{\bar{x}}^k\|^2\\
&\qquad +(1+\frac{1}{\epsilon^k\rho})(\epsilon^k\lambda^k)^2\|
 g^k\|^2\\
& \leq (1- \epsilon^k \rho) \|
x^k-\hat{\bar{x}}^k\|^2+\left((\epsilon^k\lambda^k)^2+\frac{\epsilon^k
(\lambda^k)^2}{\rho}\right) \|g^k\|^2\\
& \leq (1-\epsilon^k \rho) \|
x^k-\hat{\bar{x}}^k\|^2+\left((\epsilon^k\lambda^k)^2+\frac{\epsilon^k
(\lambda^k)^2}{\rho}\right) G^2
\end{aligned}
\end{equation}

Note that  there always exists a $\beta>0$ such that
$\mathbb{E}\left[\|V^k\|^2\right]< \beta\|x^k\|^2$ holds under
Assumption \ref{Assmption:quantization}, we can combine
(\ref{eq:consensus_first}) and (\ref{eq:Consensus_error}) to obtain
\begin{equation}\label{eq:consensus_second}
\begin{aligned}
&\mathbb{E}\left[\|x^{k+1}-\hat{\bar{x}}^{k+1}\|^2\right]
  \leq (1-\epsilon^k \rho) \mathbb{E}\left[\|
x^k-\hat{\bar{x}}^k\|^2\right]\\
&\qquad\qquad +\left((\epsilon^k\lambda^k)^2+\frac{\epsilon^k
(\lambda^k)^2}{\rho}\right)G^2
+(\epsilon^k)^2\beta\mathbb{E}\left[\| x^k \|^2\right]
\end{aligned}
\end{equation}

{\bf Step III}: Finally, combining (\ref{eq:x_bar_x_tilde_x}),
(\ref{eq:after_applying_first}), and (\ref{eq:consensus_second})
yields
\begin{equation}
\begin{aligned}
\mathbb{E}&\left[\|x^{k+1}-\hat{\bar{x}}^{k+1}\|^2+\|\bar{x}^{k+1}-x^{\ast}\|^2\right]\\&\leq
(1-\epsilon^k \rho) \mathbb{E}\left[\|
x^k-\hat{\bar{x}}^k\|^2\right]+\left((\epsilon^k\lambda^k)^2+\frac{\epsilon^k
(\lambda^k)^2}{\rho}\right)G^2 \\
&\qquad +(\epsilon^k)^2\beta\mathbb{E}\left[\| x^k
\|^2\right]+(1+(\epsilon^k)^2)
\mathbb{E}\left[\left\|\bar{x}^k-x^{\ast}\right\|^2\right] \\
&\qquad +
 \left((\epsilon^k\lambda^k)^2+(\lambda^k)^2\right)G^2 \\
&\leq (1+(\epsilon^k)^2) \mathbb{E}\left[\|
x^k-\hat{\bar{x}}^k\|^2+\left\|\bar{x}^k-x^{\ast}\right\|^2\right]\\
&\qquad +(\epsilon^k)^2\beta\mathbb{E}\left[\| x^k
\|^2\right]\\
&\qquad +\left(2(\epsilon^k\lambda^k)^2+(1+\frac{\epsilon^k
}{\rho})(\lambda^k)^2\right)G^2 \\
&\leq (1+(\epsilon^k)^2+4\beta(\epsilon^k)^2) \mathbb{E}\left[\|
x^k-\hat{\bar{x}}^k\|^2+\left\|\bar{x}^k-x^{\ast}\right\|^2\right]\\
&\qquad +\left(2(\epsilon^k\lambda^k)^2+(1+\frac{\epsilon^k
}{\rho})(\lambda^k)^2\right)G^2+2\beta(\epsilon^k)^2)\|\hat{x}^{\ast}\|^2
\end{aligned}
\end{equation}

Because the second and third terms    on the right hand side of the
above inequality are summable under the conditions in Lemma
\ref{lemma:lemma_boundedx^2}, according to Lemma
\ref{Lemma:convergence_lemma} in the Appendix, we have that
$\mathbb{E}\left[\|x^{k+1}-\hat{\bar{x}}^{k+1}\|^2+\|\bar{x}^{k+1}-x^{\ast}\|^2\right]$
will converge to a finite value. Further using
(\ref{eq:x_bar_x_tilde_x}) and the fact that $x^{\ast}$ is a finite
vector, we have that $\mathbb{E}\left[\|x^k\|^2\right]$ is always
bounded.

\subsection{\bf Proof of Lemma \ref{Lemma:convergence_of_x_i_to_bar_x}}\label{Appendix:Lemma_convergence}
  Noting that
$\mathbb{E}\left[\|x^k\|^2\right]$ is bounded from Lemma
\ref{lemma:lemma_boundedx^2}, we always have the following
inequality for some $\beta>0$ according to
(\ref{eq:consensus_second}):
\begin{equation}\label{eq:consensus_third}
\begin{aligned}
\mathbb{E}\left[\|x^{k+1}-\hat{\bar{x}}^{k+1}\|^2\right]&\leq
(1-\epsilon^k \rho) \mathbb{E}\left[\|
x^k-\hat{\bar{x}}^k\|^2\right]\\
&
 +\left((\epsilon^k\lambda^k)^2+\frac{\epsilon^k(\lambda^k)^2 }{\rho}\right)G^2
+(\epsilon^k)^2\beta \Omega
\end{aligned}
\end{equation}
where $\Omega$ is some constant representing an upper bound of
$\mathbb{E}\left[\|x^k\|^2\right]$. Then the lemma can be directly
obtained by applying Lemma \ref{Lemma:decay_rates} in Appendix A.

\bibliographystyle{unsrt}
\bibliography{reference1}

\end{document}